\documentclass[aop,seceqn,MSNbibl,citesort,dvips]{arximspdf}
\usepackage{mathbh}

\doi{10.1214/10-AOP575}
\volume{39}
\issue{3}
\pubyear{2011}
\firstpage{985}
\lastpage{1026}

\makeatletter

\newcommand{\Prr}{\operatorname{Pr}}

\newcommand{\R}{\mathbb{R}}
\newcommand{\sgn}{\operatorname{sgn}}
\newcommand{\Z}{\mathbb{Z}}
\newcommand{\E}{\mathbb{E}}
\newcommand{\eps}{\varepsilon}
\newcommand{\ka}{\kappa}
\newcommand{\VC}{\operatorname{VC}}

\newtheorem{Theorem}{Theorem}[section]
\newtheorem{Lemma}[Theorem]{Lemma}
\newproclaim{Definition}[Theorem]{Definition}
\newtheorem{Corollary}[Theorem]{Corollary}
\newproclaim{Remark}[Theorem]{Remark}

\newcommand{\IND}{\mathbh{1}}

\makeatother

\begin{document}
\begin{frontmatter}

\title{Discrepancy, chaining and subgaussian processes}
\runtitle{Discrepancy, chaining and subgaussian processes}

\begin{aug}
\author[A]{\fnms{Shahar} \snm{Mendelson}\corref{}\ead
[label=e1]{shahar.mendelson@anu.edu.au}\ead
[label=e2]{shahar@tx.technion.ac.il}}
\runauthor{S. Mendelson}
\affiliation{Israel Institute of Technology and The Australian
National University}
\address[A]{Department of Mathematics\\
Technion, Israel Institute of Technology\\
Haifa 32000\\
Israel\\
\printead{e1}\\
and\\
Centre for Mathematics and its Applications\\
Institute for Advanced Studies\\
The Australian National University\\
Canberra, ACT 0200\\
Australia\\
\printead{e2}} 
\end{aug}

\received{\smonth{10} \syear{2008}}
\revised{\smonth{6} \syear{2010}}

%
\begin{abstract}
We show that for a typical coordinate projection of a subgaussian class
of functions, the infimum over signs $\inf_{(\eps_i)} {\sup_{f \in
F}} |{\sum_{i=1}^k \eps_i} f(X_i)|$ is asymptotically smaller than the
expectation over signs
as a function of the dimension $k$, if the canonical Gaussian process
indexed by $F$ is continuous. To that end, we establish a bound on
the discrepancy of an arbitrary subset of $\R^k$ using properties of
the canonical Gaussian process the set indexes, and then obtain
quantitative structural information on a typical coordinate
projection of a subgaussian class.
\end{abstract}

%
\begin{keyword}[class=AMS]
\kwd{60C05}
\kwd{60G15}
\kwd{60D05}.
\end{keyword}
\begin{keyword}
\kwd{Discrepancy}
\kwd{generic chaining}.
\end{keyword}

\end{frontmatter}

\section{Introduction}\label{intro}
The geometric structure of a random coordinate projection of a class of
functions plays an important role in Empirical Processes theory,
where it is used to determine whether the uniform law of large
numbers or the uniform central limit theorem is satisfied by the given class.
Indeed,
if $F$ is a class of functions on a probability space
$(\Omega,\mu)$, and if $\sigma=(X_1,\ldots,X_k)$ is an independent
sample distributed according to $\mu^k$, then the ``complexity'' of
the random set
\[
P_\sigma F=\{(f(X_1),\ldots,f(X_k)) \dvtx f \in F\} \subset\R^k
\]
is the key parameter in addressing both these questions. In this
context, if $(\eps_i)_{i=1}^k$ are independent, symmetric,
$\{-1,1\}$-valued random variables, then
the complexity is governed by the expectation
of the supremum of the Bernoulli process indexed by $P_\sigma F$,
defined by
%
\begin{equation} \label{eq:intro-Bernoulli}
\E_\eps\sup_{f \in F} \Biggl|\sum_{i=1}^k \eps_i f(X_i)
\Biggr|=\E_\eps\sup_{v \in P_\sigma F} \Biggl|\sum_{i=1}^k \eps_i
v_i \Biggr|,
\end{equation}
and in particular, on the way
this expectation grows as a function of $k$ for a typical sample
of cardinality $k$ (see, e.g.,
\cite{Dud-book,VW,LT} and references therein).

The structure of such coordinate projections is central to questions in
Asymptotic
Geometric Analysis as well. For example, let $K \subset\R^d$ be a convex,
symmetric set (i.e., if $x \in K$ then $-x \in K$) and put $F=\{\langle
x,\cdot \rangle \dvtx x \in K\}$ to be the
class of linear functionals indexed by $K$. If $\mu$ is a measure on
$\R^d$, then $P_\sigma F$ is the set $\Gamma K$, where $\Gamma$ is
the random
operator $\Gamma= \sum_{i=1}^k \langle X_i,\cdot \rangle e_i$. Fundamental
questions on the geometry of convex, symmetric sets, such as
Dvoretzky's theorem and low-$M^*$ estimates have been answered by
obtaining accurate, quantitative information on the structure of such coordinate
projections, and by
using very similar complexity parameters to (\ref{eq:intro-Bernoulli})
(e.g., \cite{MS,LT}).

For both these reasons, a lot of effort has been invested in
understanding various notions of complexity for
a typical coordinate projection of a class of functions.
A well studied direction is to obtain
quantitative estimates on the way in
which (\ref{eq:intro-Bernoulli}) is related to two other complexity
parameters, the combinatorial dimension and covering numbers.

Roughly speaking, the combinatorial dimension of $T \subset\R^k$ at
scale $\eps$, denoted by $\operatorname{VC}(T,\eps)$, is the largest dimension of a
coordinate projection of $T$ that contains a ``cube'' of scale
$\eps$ (see Definition \ref{def:comb-dim} for an exact formulation).
If $(T,d)$ is a metric space then the covering number at scale
$\eps$, which we denote by $N(\eps,T,d)$, is the smallest cardinality
of a
subset $\{y_1,\ldots,y_m\} \subset T$ such that for every $t \in T$,
there is some $y_i$ for which $d(t,y_i) < \eps$.

Connections between (\ref{eq:intro-Bernoulli}) and the
combinatorial dimension or the covering numbers of $P_\sigma F$ are
rather well
understood. For example, a straightforward chaining argument (see,
e.g., \cite{VW})
shows that for every $T \subset\R^k$,
%
\begin{equation} \label{eq:intro-entropy-integral}
\E_\eps\sup_{t \in T} \Biggl|\sum_{i=1}^k \eps_i t_i\Biggr| \leq
c\int_0^{\operatorname{diam}(T)} \sqrt{\log N(\eps,T,\ell_2^k)}\,d\eps,
\end{equation}
where $\ell_2^k$ is the Euclidean metric on $\R^k$,
$\operatorname{diam}(T)$ is the diameter with respect to the same
metric and $c$ is an absolute constant, independent of the dimension
$k$ and of the set $T$. Some of the other relations between these
parameters are far more involved. First, controlling the $L_2$ covering
numbers using the combinatorial dimension was resolved in
\cite{MenVer}, where it was shown that if $T$ is a subset of the unit
cube $B_\infty^k$ and $\mu$ is any probability measure on
$\{1,\ldots,k\}$, then for every $0<\eps<1$,
\[
N(\eps,T,L_2(\mu)) \leq\biggl(\frac{5}{\eps}\biggr)^{K \cdot
\operatorname{VC}(T,c\eps)},
\]
where $K$ and $c$ are absolute constants. Also, the solution of the
sign embedding of $\ell_1^k$ problem (see \cite{MenVer}) implies that
if $T \subset B_\infty^k$ and ${\E\sup_{t \in T}} |{\sum_{i=1}^k \eps_i
t_i} | \geq\delta k$, then $\operatorname{VC}(T,c_1\delta) \geq c_2\delta^2 k$. In
other words, under a normalization condition ($T \subset B_\infty^k$),
the only reason that ${\E\sup_{t \in T}} |{\sum_{i=1}^k \eps_i t_i }|$ is
almost extremal is that $T$ contains a high-dimensional cubic
structure.

In this article, we study a related geometric parameter---the
discrepancy of a typical coordinate projection. Discrepancy was
introduced in a combinatorial context (see the book \cite{Matbook} for
an extensive survey on this topic) and is defined as follows.
\begin{Definition} \label{def:disc}
If $T \subset\R^k$, then the discrepancy of $T$ is
\[
\operatorname{disc}(T)=\inf_{(\eps_i)_{i=1}^k} \sup_{t \in T} \Biggl|
\sum_{i=1}^k \eps_i t_i \Biggr|,
\]
and the infimum is taken with respect to all signs
$(\eps_i)_{i=1}^k \in\{-1,1\}^k$.

We denote by $\operatorname{Hdisc}(T)$ the {hereditary discrepancy} of $T$,
given by
\[
\sup_{I \subset\{1,\ldots,k\}} \operatorname{disc}(P_I T),
\]
where $P_I T=\{(t_i)_{i \in I} \dvtx t \in T\}$ is the
coordinate projection of $T$ onto $I$.
\end{Definition}

Observe that if $\operatorname{absconv}(T)$ is the convex hull of $T \cup-T$,
then $\operatorname{disc}(T)=\operatorname{disc}(\operatorname{absconv}(T))$. Hence, from the
geometric viewpoint, the discrepancy of $T$ is proportional with a
constant $\sqrt{k}$
to the
minimal width of $\operatorname{absconv}(T)$ in a direction of a vertex of the
combinatorial cube $\{-1,1\}^k$. The hereditary discrepancy is governed
by a
similar minimal width, but of the ``worst''
coordinate projection of $\operatorname{absconv}(T)$.

Our goal here is to study the discrepancy using
the covering numbers and the combinatorial dimension of $T$, but
we will focus on sets $T$ that
are random coordinate projections of a class of function $F$, which
gives them
more structure.
A~natural question in this context is to identify
conditions on $F$ under which there is a gap between $\operatorname{disc}(P_\sigma F)$
and ${\E_\eps\sup_{v \in P_\sigma F}} |{\sum_{i=1}^k \eps_i v_i}|$ for
a typical
coordinate projection of $F$, as a function of the sample size $k$.
To that end, we
will develop dimension dependent bounds on the discrepancy of $P_\sigma
F$ (and
in particular, bounds that are not asymptotic).

Note that the metric structure of $T \subset\R^k$ is not enough to
determine its
discrepancy. Indeed, if $e_1=(1,0,\ldots,0) \in\R^k$ and
$T_1=\{0,e_1\}$ then $\operatorname{disc}(T_1)=1$. On the other hand, if
$T_2=\{0,\sum_{i=1}^k e_i/\sqrt{k} \}$, which is linearly isometric
to $T_1$,
then $\operatorname{disc}(T_2) \leq1/\sqrt{k}$. The reason for the large gap
in the discrepancy between the two isometric sets is that
$T_2$ consists of a vector that is ``well
spread'' while $T_1$ consists of a ``peaky'' vector with respect to the
underlying coordinate structure. In that sense,
$T_2$ is in a much better position than $T_1$. Note that in this example,
${\E\sup_{t \in T_1} }|{\sum_{i=1}^k \eps_it_i}| =
{\E\sup_{t \in T_2}} |{\sum_{i=1}^k \eps_it_i}|=1$---and for the set $T_2$,
which is in a ``good position'' there is gap between the expectation for signs
and the discrepancy.

We will show that this is a general phenomenon:
it is well known that ${\E_\eps\sup_{t \in T} }|{\sum_{i=1}^k \eps_i
t_i}|$ is determined
by the Euclidean metric structure of $T$ (up to a logarithmic factor in the
dimension), and therefore, it is almost invariant under a linear isometry
(i.e., a change in the coordinate structure).
Thus, the expectation almost does not change when applying
an isometry or a good isomorphism of $\ell_2^k$. As we will explain
here, the situation
with the discrepancy is rather different and the position of the set
matters a great deal.
Since the sets $T$
that we will be interested in are not arbitrary but
have some structure---as random coordinate projections of well
behaved function classes, they will be much closer in nature to
$T_2$ than to $T_1$.

Our main result is that if
the canonical Gaussian process indexed by $F \subset L_2(\mu)$ is continuous
and if the class satisfies a subgaussian condition [i.e., if the
$\psi_2(\mu)$ norm is equivalent to the
$L_2(\mu)$ norm on $F$, see Definition \ref{def:sub-gauss-set}], then
a typical
coordinate projection of $F$ behaves as a set of vectors
in a ``general position.'' As such, and just like the set $T_2$, a
typical coordinate projection exhibits certain shrinking
properties that will be explained in Section \ref{sec:shrink}, and
which causes the discrepancy
of such a set to be much smaller than the average over signs.

{
\setcounter{Theorem}{0}
\renewcommand{\theTheorem}{\Alph{Theorem}}
\begin{Theorem}\label{theoA}
Let $F \subset L_2(\mu)$ be a class of mean zero functions. Assume further
that the canonical Gaussian process indexed by $F$ is continuous and that
the $L_2(\mu)$ and $\psi_2(\mu)$ norms are equivalent on $F$.
Then $\operatorname{Hdisc}(P_\sigma F)/\sqrt{k} \to0$ in probability.
\end{Theorem}

To put Theorem \ref{theoA} in the right perspective, observe that
if the $\psi_2$ and $L_2$ norms are equivalent on a class of mean zero
functions $F$,
then for every integer $k$ there is a subset of $\Omega^k$ of
probability at least $c$ on which
%
\begin{equation} \label{eq:lower}
\E_\eps\sup_{f \in F}
\Biggl| \sum_{i=1}^k \eps_i f(X_i) \Biggr| \geq c_1\sqrt{k}\sigma_F,
\end{equation}
where $c$ depends only on the equivalence constant between the $\psi
_2$ and $L_2$ norms on $F$,
$c_1$ is an absolute constant and $\sigma_F =\sup_{f \in F} (\E
f^2)^{1/2}$. Hence, there is a true gap between the discrepancy
and the
mean of a typical coordinate projection.

Although the formulation of Theorem \ref{theoA} is asymptotic, the
result itself is quantitative in nature, as a function of the dimension
of the coordinate projection. The proof of Theorem \ref{theoA} is, in
fact, a dimension dependent estimate on the sequences
$(\alpha_{k,\delta})_{k=1}^\infty$, for which, with probability at
least $1-\delta$, $\operatorname{Hdisc}(P_\sigma F)
\leq\sqrt{k}a_{k,\delta}$. We will show that the sequences
$(\alpha_{k,\delta})_{k=1}^\infty$ are given using metric parameters
that measure the continuity of the Gaussian process indexed by
$F$---Talagrand's $\gamma_{2,s}$ functionals \cite{Tal-book}. The
$\gamma_{2,s}$ functionals will be defined in Section~\ref{sec:pre},
but for now let us mention that under mild measurability assumptions on
the class, the canonical Gaussian process indexed by $F$ is continuous
if and only if $\lim_{s \to\infty} \gamma_{2,s}(F,L_2(\mu)) = 0$.

We will prove that for every $0<\rho<1/2$ and $0<\delta<1$ there are
constants $c$ and $C$ that depend on $\rho$, $\delta$ and on the
equivalence constant between
the $\psi_2(\mu)$ and $L_2(\mu)$ norms on $F$, such that for every $k$,
%
\begin{equation} \label{eq:est-a_k}\hspace*{28pt}
\alpha_{k,\delta} \leq C \sup_{1 \leq n \leq k} \sqrt{\frac{n}{k}}
\biggl(\gamma_{2,\log_2 \log_2 cn} (F,L_2(\mu)) \cdot\sqrt{\log
(ek/n)} +
D\frac{\log{k}}{n^{{1/2}-\rho}}\biggr),
\end{equation}
where $D=\operatorname{diam}(F,L_2(\mu))$ is the diameter of $F$ with respect
to the $L_2(\mu)$ norm. And, in particular, under the assumptions of
Theorem \ref{theoA},
for every $0<\delta<1$, $\lim_{k \to\infty} \alpha_{k,\delta}=0$.
Moreover, the proof of Theorem \ref{theoA} actually shows
that for every~$k$, with $\mu^k$-probability of at least $1-\delta$,
\[
\operatorname{disc}(P_\sigma F) \leq C \bigl(
\sqrt{k}\gamma_{2,\log_2 \log_2 ck} (F,L_2(\mu))+k^\rho D\bigr),
\]
where $C$ and $c$ depend on $\rho$, $\delta$ and the equivalence constant
between the $L_2(\mu)$ and $\psi_2(\mu)$ norms on $F$.

The proof of Theorem \ref{theoA} is based on two ingredients. The first is a
new estimate on the discrepancy of an arbitrary set $T \subset
\R^k$. It is a combination of the entropy method, which is often
used to control the combinatorial discrepancy (see, e.g.,
\cite{Spen,AS,Matbook}), and Talagrand's generic chaining mechanism
\cite{Tal-book}, which was introduced to establish the connection
between the $\gamma_{2,s}$ functionals and the continuity of
Gaussian processes. The combination of these two methods will be
explained in Section \ref{sec:disc}. It allows one to find
a good choice of signs for roughly $k/2$ coordinates, while the
error incurred by considering the sum taken only on these
coordinates is determined by the $\gamma_{2,s}$ functional for
$s \sim\log_2 \log_2 k$.
Repeating this argument, one obtains a bound on the
discrepancy of $T$ in terms of a sum of $\gamma_{2,s}$ functionals of
coordinate
projections of $T$ and for values $s$ that depend on the dimension of
each projection,
and those dimensions decrease quickly.

The second component required for the proof of Theorem \ref{theoA} is that the
sets we
are interested in are not general. We will obtain a
structural result on a typical $P_\sigma F$ that allows us to bound the
$\gamma_{2,s}$
functionals of its coordinate projections using the $L_2(\mu)$
structure of $F$.

Indeed, we will show that if the $L_2(\mu)$ and $\psi_2(\mu)$ norms
are equivalent on $F$ then
a typical coordinate projection $P_\sigma F$ has a rather regular
structure---it is
a subset of a Minkowski sum of two sets. The first one is
small, with a bounded diameter in the weak $\ell_2$ space
$\ell_{2,\infty}^k$. The other set is a subset of $P_\sigma F$
itself and can be viewed as a set of vectors in a ``general position.''
In particular, further coordinate projections of the latter set
shrink distances between any two of its elements.

The structural result we obtain is of independent interest
and can be used to derive information on the geometry of
convex sets. For example,
consider a symmetric probability measure $\mu$ on $\R^n$.
We say that $\mu$
isotropic and $L$-subgaussian if a random vector $X$ distributed
according to $\mu$ satisfies that for every $x \in\R^n$,
\[
\E|\langle X,x \rangle|^2=|x|^2  \quad\mbox{and}\quad \|\langle X,x
\rangle\|_{\psi_2} \leq
L|x|.
\]
Simple examples of isotropic, $L$-subgaussian measures on $\R^n$ are
the Gaussian measure and the uniform measure on the vertices of the
cube $\{-1,1\}^n$, where in both cases $L$ can be taken to be an absolute
constant, independent of the dimension.

Let $(X_i)_{i=1}^k$ be independent random vectors,
distributed according to $\mu$ and consider the random operator
$\Gamma\dvtx\R^n\to\R^k$ defined by $\Gamma=\sum_{i=1}^k
\langle X_i,\cdot \rangle e_i$.
\begin{Corollary}\label{corolB}
For any $L>0$ there are constants $c_1,c_2$ and $c_3$
that depend only on $L$, for which the following holds.
Let $T \subset\R^n$ and set
$V=k^{-1/2}\Gamma T$. Then, for every $u>c_1$, with probability at
least $1-2\exp(-c_2u)$, for every $I \subset\{1,\ldots,k\}$,
\[
\E_g \sup_{v \in V} \biggl| \sum_{i \in I} g_i v_i \biggr|
\leq c_3 u \sqrt{\frac{|I|}{k}\log\biggl(\frac{ek}{|I|}\biggr)}
\E_g \sup_{t\in T} \Biggl|\sum_{i=1}^n g_i t_i\Biggr|,
\]
where $(g_i)$ are independent, standard Gaussian random variables,
and both expectations are taken with respect to those variables.
\end{Corollary}
}

Corollary \ref{corolB} shows that the random operator $\Gamma$ maps an arbitrary $T$
to a set of vectors in a ``general position'' in a strong sense, since
it implies that for most vectors in $V$, mutual distances are shrunk by
any further coordinate projection. Let us note that we will prove a
stronger result than Corollary \ref{corolB}, namely that the $\gamma_{2,s}$ functionals
associated with $V$ display this type of shrinking phenomenon.

The final result we present has to do with the reverse direction of
Theorem~\ref{theoA}.
Assume that $H \subset L_2(\mu)$ is a convex, symmetric set,
which satisfies that the canonical Gaussian process $\{G_h \dvtx h \in H\}$
is bounded
and that the $L_2(\mu)$ and $\psi_2(\mu)$ norms are equivalent on $H$.
We will
show that if the logarithm of the $L_2(\mu)$ covering numbers of $H$
grows like $1/\eps^2$ then for a typical sample
$\sigma=(X_1,\ldots,X_k)$ selected according to $\mu^k$,
\[
\operatorname{VC}\bigl(P_\sigma H, c_1/\sqrt{k}\bigr) \geq c_2k.
\]
It is standard to verify (see Lemma \ref{lemma:disc-vs-VC}) that if
$T \subset\R^k$, then
\[
\operatorname{Hdisc}(T) \geq\sup_{\delta>0} \delta\operatorname{VC}(\operatorname{absconv}(T),\delta).
\]

Therefore, if $F$ is a class of mean-zero functions and $H=\operatorname{absconv}(F)$ satisfies
the above, then $\operatorname{Hdisc}(P_\sigma F) \geq c\sqrt{k}$,
complementing the upper bound established in Theorem \ref{theoA}.

Although this is not exactly the reverse direction of Theorem \ref{theoA}, it
is very close to it. Indeed, if $F \subset L_2(\mu)$ indexes a
bounded Gaussian process then so does $H=\operatorname{absconv}(F)$, and the
logarithm of the covering numbers
of $H$ cannot grow faster than $O(1/\eps^2)$. On the other hand, if
the log-covering numbers grow a little slower, even by a suitable
logarithmic factor, then
$\gamma_{2,s}(F,L_2(\mu)) \to0$. In fact, this is as close as one
can get to a covering numbers characterization of the fact that
$\gamma_{2,s}(F,L_2) \to0$ (see, e.g., \cite{Dud-book}).

This result not only shows that $\operatorname{Hdisc}(P_\sigma F)$ is large if
the Gaussian
process $F$ indexes is bounded but not continuous, it also shows why.
Under a boundedness
assumption on the Gaussian process [which implies that $\operatorname{Hdisc}(P_\sigma F)/\sqrt{k}$
is bounded], the reason the hereditary discrepancy of $P_\sigma F$ is
extremal is because
a typical coordinate projection of $\operatorname{absconv}(F)$ contains a high
dimensional, large cubic structure, and that forces the hereditary
discrepancy to be large.
The proof of this result, which is presented in Section \ref{sec:upp-low},
is based on the observation that if $F$ is convex and symmetric then
the richness of $F$ at scale
$\sim1/\sqrt{k}$ is exhibited by the existence of cubes of scale
$\sim1/\sqrt{k}$ and of dimension $\sim k$ in a typical coordinate
projection of $F$ of dimension $k$. It thus should be viewed as a
``small scale'' version of the Sign Embedding theorem which was
mentioned above.

Unfortunately, the optimal estimate in the Sign Embedding theorem
cannot be used directly in our
case, firstly because $P_\sigma F$ is unlikely to be a subset of
$B_\infty^k$, and secondly, because a typical coordinate projection of
$F$ satisfies that
\[
\E\sup_{v \in P_\sigma F} \Biggl|\sum_{i=1}^k \eps_i v_i \Biggr|
\sim\sqrt{k}.
\]
Hence, the optimal estimate in the Sign Embedding theorem has to be
used for
$\delta\sim1/\sqrt{k}$, and that only ensures that $P_\sigma F$
contains a cube of scale $\sim1/\sqrt{k}$ and of constant dimension,
which is far from what we need.

The proof of the existence of a cube in $P_\sigma F$ is based on two
localization arguments, one
with respect to the $L_2$ norm and the other with respect to the
$L_\infty$ norm. The first localization shows that if the $L_2(\mu)$
covering number of $F$ at scale $\sim1/\sqrt{k}$ is of the order
of $\exp(c_1k)$ then the richness of a typical coordinate
projection of $F$ of dimension $\sim k$ originates from the set
%
\begin{equation} \label{eq:intersect}
F_1=F \cap\frac{c_2}{\sqrt{k}}B(L_2(\mu)),
\end{equation}
that is, functions in $F$ of $L_2(\mu)$ norm at most
$O(1/\sqrt{k})$. In the second localization, one shows that the
complexity of a typical coordinate projection actually comes from a
further pointwise truncation of the functions in $F$, and
$B(L_2(\mu))$ in (\ref{eq:intersect}) can essentially be replaced by
$B(L_\infty(\mu))$---the unit ball in $L_\infty(\mu)$.

This article is organized as follows. In Section \ref{sec:pre}, we
present further preliminaries, most of them concerning subgaussian
variables and the $\gamma_{2,s}$ functionals. In Section
\ref{sec:disc}, we develop bounds on the discrepancy of an arbitrary
subset of $\R^n$. Section \ref{sec:shrink} is devoted to the proof of
the structural results on coordinate projections of subgaussian
processes and its corollaries, including Corollary \ref{corolB}.
Theorem~\ref{theoA} is proved in Section \ref{sec:main} and its
converse and the resulting lower bound on the hereditary discrepancy of
a typical coordinate projection is proved in Section \ref{sec:upp-low}.

\section{Preliminaries} \label{sec:pre}
Throughout, absolute constants (i.e., fixed, positive numbers) will
be denoted by $C,c,c_1$ etc. Their values may change from line to
line. We denote by $C(a),c(a)$ constants that depend only on the
parameter $a$ and we set $\ka_1,\ka_2,\ldots$ to be constants that will
remain fixed throughout the article. By $a \sim b$, we mean that
there are constants $c$ and $C$ such that $ca \leq b \leq Ca$, and we
write $b \lesssim a$ if $b \leq Ca$.

We will consider a single, fixed Euclidean structure on all
finite-dimensional spaces $\R^n$ and denote the corresponding Euclidean
norms by \mbox{$| \cdot|$} without specifying the dimension. With a minor
abuse of notation, the cardinality of a set and the absolute value are
denoted in the same way.

If $E$ is a normed space, let $B(E)$ be its unit ball, and for
$E=\ell_p^n$ we denote the unit ball by $B_p^n$. If
$\sigma=(X_1,\ldots,X_k) \in\Omega^k$ let $\mu_k=k^{-1}\sum_{i=1}^k
\delta_{X_i}$ be the empirical measure supported on $\sigma$, set
$L_2^k$ to be the corresponding $L_2$ space, and for $I \subset\{
1,\ldots,k\}$
let $\ell_2^I$ be the
coordinate subspace of $\ell_2^k$ spanned by
$(e_i)_{i \in I}$.

The situation we will study here is as follows.
Let $F$ be a class of real valued functions on a
probability space $(\Omega,\mu)$, let $X_1,\ldots,X_k$ be independent
random variables distributed according to $\mu$ and set
$\sigma=(X_1,\ldots,X_k)$. Let $P_\sigma F=\{(f(X_i))_{i=1}^k \dvtx f \in
F\} \subset\R^k$ be the coordinate projection of $F$ defined by
$\sigma$
and for every $I \subset\{1,\ldots,k\}$ let
$P_I^\sigma F = \{ (f(X_i))_{i \in I} \dvtx f \in F
\} \subset\R^{|I|}$ be the coordinate projection of $F$ on the
subset of
coordinates $(X_i)_{i \in I}$. Sometimes, for the sake of simplicity,
we shall
omit the superscript $\sigma$.

\subsection{Subgaussian processes}
Here, we will describe properties of sums of independent random
variables that have quickly decaying tails.
\begin{Definition} \label{def:psi}
Let $f$ be a functions defined on a probability space $(\Omega,\mu)$.
For $1 \leq\alpha\leq2$, define the $\alpha$-Orlicz norm by
\[
\|f\|_{\psi_\alpha} =
\inf\biggl\{C>0\dvtx \E\exp\biggl(\frac{|f|^\alpha}{C^\alpha}\biggr)
\leq2 \biggr\}.
\]
\end{Definition}

For basic facts regarding Orlicz norms, we refer the reader to
\cite{dPG,VW}.

It is well known that a random
variable has a bounded $\psi_\alpha$ norm for $1 \leq\alpha\leq2$
if and only if it has a well behaved tail; that is, there is an absolute
constant $c$ such that for every $f \in L_{\psi_\alpha}$ and every $t
\geq1$,
\[
\Prr (|f| \geq t) \leq
2\exp(-ct^\alpha/\|f\|_{\psi_\alpha}^\alpha).
\]
Conversely, there is an absolute constant $c_1$ such that
if $f$ displays a tail behavior dominated by
$\exp(-t^\alpha/K^\alpha)$ for $1 \leq\alpha\leq2$ then
$\|f\|_{\psi_\alpha} \leq c_1K$.

There are several basic properties of sums of independent random
variables we require. The proofs of these facts can be found, for
example, in \cite{LT,VW,dPG}.

Note that if $f$ has a subexponential tail
then its empirical means concentrate around its true mean, with
a tail behavior that is a mixture of subgaussian and subexponential.
Indeed, the following result is a version of Bernstein's inequality
and shows just that.
\begin{Theorem} \label{thm:Bern}
There exists an absolute constant $c$ for which the following holds.
Let $(\Omega,\mu)$ be a probability space and set $f\dvtx\Omega\to\R$
to be a function with a bounded $\psi_1$ norm. If $X_1,\ldots,X_k$ are
independent and distributed according to $\mu$ then for every $t>0$,
\[
\Prr \Biggl( \Biggl|\frac{1}{k}\sum_{i=1}^k f(X_i)-\E f \Biggr| \geq
t\|f\|_{\psi_1} \Biggr) \leq2\exp(-ck\min\{t^2,t\}).
\]
\end{Theorem}

If a function has mean zero and a bounded $\psi_2$ norm, one can
obtain a purely subgaussian tail.
\begin{Lemma} \label{Lemma:psi_2}
There exists an absolute constant $c$ for which the following holds.
Let $Y_1,\ldots,Y_k$ be independent random variables of mean zero.
Then, for every $a_1,\ldots,a_k \in\R$,
\[
\Biggl\|\sum_{i=1}^k a_i Y_i \Biggr\|_{\psi_2} \leq
c\Biggl(\sum_{i=1}^k a_i^2 \|Y_i\|_{\psi_2}^2 \Biggr)^{1/2}.
\]
In particular, if $(X_i)_{i=1}^k$ are independent random variables
distributed according to $\mu$ and $f$ has zero mean, then for every
$t \geq1$,
\[
\Prr \Biggl(\Biggl|\sum_{i=1}^k f(X_i) \Biggr| \geq t
k^{1/2}\|f\|_{\psi_2} \Biggr) \leq2\exp(-c_1t^2),
\]
where $c_1$ is an absolute constant.
\end{Lemma}

In what follows, we will assume that the $\psi_2$ and $L_2$ norms are
equivalent on $F$
in the following sense.
\begin{Definition} \label{def:sub-gauss-set}
A set $F \subset L_2(\mu)$ is $L$-subgaussian if $\|f\|_{\psi_2} \leq
L\|f\|_{L_2}$ and $\|f-g\|_{\psi_2} \leq L \|f-g\|_{L_2}$ for every
$f,g \in F$.
\end{Definition}

Next, let us turn to the definition of the $\gamma_{2,s}$ functionals
\cite{Tal-book}.
Let $(T,d)$ be a metric space. An \textit{admissible sequence} of $T$ is a sequence of subsets of $T$,
$\{T_s\}_{s=0}^\infty$, such that $|T_0|=1$ and for every $s \geq
1$, $|T_s| \leq2^{2^s}$.
\begin{Definition} \label{def:gamma-2}
For a metric space $(T,d)$ and an integer $s_0\geq0$, let
\[
\gamma_{2,s_0}(T,d) =\inf\sup_{t \in T} \sum_{s=s_0}^\infty2^{s/2}d(t,T_s),
\]
where the infimum is taken with respect to all admissible sequences of $T$.
Set $\gamma_2(T,d) =\gamma_{2,0}(T,d)$.
\end{Definition}

Let $\pi_s\dvtx T \to T_s$
be a metric projection function onto $T_s$, that is, $\pi_s(t)$
is a nearest point to $t$ in $T_s$ with respect to the metric $d$.
It is easy to verify that for every admissible sequence, every
$t \in T$, and any $s_0 \geq0$,
\[
\sum_{s=s_0}^\infty2^{s/2}d(\pi_{s+1}(t),\pi_{s}(t)) \leq\bigl( 1 +
1/\sqrt2\bigr) \sum_{s=s_0}^\infty2^{s/2}d(t,T_s)
\]
and that the diameter of $T$ satisfies $\operatorname{diam}(T,d) \leq
2\gamma_2(T,d)$. Moreover, it is clear that the $\gamma_{2,s}$
functionals are decreasing in $s$ and are subadditive in $T$ in the
following sense. Let $X$ be a normed space and consider two sets $A,B
\subset X$. If
$A+B=\{a+b \dvtx a \in A, b\in B\}$ is the Minkowski sum of $A$ and
$B$, then for every integer $s$,
\[
\gamma_{2,s+1}(A + B,d) \leq\gamma_{2,s}(A,d)+\gamma_{2,s}(B,d).
\]
There is a close connection between the $\gamma_{2,s}$ functionals with
respect to $L_2$ norms and properties of Gaussian processes (see
\cite{Dud-book,Tal-book} for expositions on these connections).
Indeed, let $\{G_t\dvtx t \in T\}$ be a centered Gaussian process
indexed by a set $T$ and for every $s, t \in T$ define a metric on
$T$ by $d^2 (s,t) = \E|G_s-G_t|^2$. One can show that
under mild measurability assumptions on $T$,
\[
c_1 \gamma_2(T,d) \leq\E\sup_{t \in T} G_t \leq c_2 \gamma_2(T,d),
\]
where $c_1$ and $c_2$ are absolute constants. The upper bound is
due to Fernique \cite{F} and the lower bound is Talagrand's
Majorizing Measures theorem \cite{Tal87}. The proof of both parts
can be found in \cite{Tal-book}. Thus, the $\gamma_2$ functional is
finite if and only if the Gaussian process indexed by $T$ is
bounded.

Note that if $T \subset\R^n$ and $G_t = \sum_{i=1}^n g_i t_i$ then
$d(u, t) = |u-t|$ and therefore
%
\begin{equation}
\label{maj_meas}
c_1 \gamma_2(T,| \cdot|) \leq\E\sup_{t \in T} \sum_{i=1}^n
g_it_i \leq c_2 \gamma_2(T,|\cdot|).
\end{equation}

Just like $\gamma_2(T,L_2(\mu))$ determines the supremum of the
canonical Gaussian process indexed by $T \subset L_2(\mu)$ (which we will
always assume to satisfy the necessary measurability assumptions), the
continuity of that process is determined by properties of the
sequence $\gamma_{2,s}$.
\begin{Definition} \label{def:pregaussian}
Let $F \subset L_2(\mu)$ be a class of mean zero functions. Set
$\{G_f \dvtx f \in F\}$ to be the centered Gaussian process indexed by
$F$ with a covariance structure endowed by $L_2(\mu)$; that is, for
every $f,g \in F$, $\E G_f G_g =\langle f,g \rangle_{L_2(\mu)}$. We
say that
$F$ is $\mu$-pregaussian if it has a version with all sample
functions bounded and uniformly continuous with respect to the
$L_2(\mu)$ metric.
\end{Definition}
\begin{Theorem}[\cite{Tal87,Tal-book}] \label{thm:cont-gauss}
Let $\{G_t \dvtx t \in T\}$ be a centered Gaussian process and endow $T$
with the $L_2$ metric given by the covariance structure of the
process, as above. Under measurability assumptions, the following are
equivalent:
\begin{enumerate}
\item The map $t \to G_t(\omega)$ is uniformly continuous on $T$
with probability $1$.
\item $\lim_{\delta\to0} \E\sup_{d(u,t) \leq\delta} |G_u-G_t|
= 0$.
\item There exists an admissible sequence of $T$ such that
\[
\lim_{s_0 \to\infty} \sup_{t \in T} \sum_{s=s_0}^\infty
2^{s/2}d(t,\pi_s(t))=0.
\]
\end{enumerate}
\end{Theorem}

In other words, $T$ is pregaussian if and only if $\lim_{s \to\infty
} \gamma_{2,s}(T,L_2)=0$.
\begin{Remark}
Theorem \ref{thm:cont-gauss} is not proved in \cite{Tal-book}
but only stated there, and its formulation in \cite{Tal87} was done
using the notion of majorizing measures rather than with the
$\gamma_{2,s}$ functionals. Since the proof of the continuity
theorem follows from an application of the Majorizing Measures
theorem and since the latter is proved in \cite{Tal-book} using the
language of the $\gamma_2$-functional, it is not difficult to
convert the proof of the continuity theorem from \cite{Tal87} and
obtain Theorem \ref{thm:cont-gauss}. Moreover, as shown in
\cite{Tal87}, there is a quantitative connection between the modulus
of continuity of $\{G_t \dvtx t \in T\}$ and the sequence
$(\gamma_{2,s}(T,L_2))_{s=0}^\infty$.
Since we will not use this
quantitative estimate here, we will not formulate it.
\end{Remark}

Finally, let us define the covering and packing numbers of a metric space.
\begin{Definition} \label{def:cov-and-ent}
Let $(T,d)$ be a metric space. The covering number of $T$ at scale
$\eps>0$ with respect to the metric $d$ is the smallest number of
open balls of radius $\eps$ needed to cover $T$, and is denoted by
$N(\eps,T,d)$.

We set $e_k(T,d)=\inf\{\eps\dvtx N(\eps,T,d) \leq2^k\}$.
$(e_k)_{k=0}^\infty$ are called the entropy numbers of $T$.

A set $A \subset T$ is called $\eps$-separated if the distance between
any two of its elements is at least $\eps$.
We denote by $D(\eps,T,d)$ the cardinality of a maximal $\eps
$-separated subset of $T$.
\end{Definition}

It is standard to verify that for every $\eps>0$, $N(\eps,T,d) \leq
D(\eps,T,d) \leq N(\eps/2$, $T,d)$,
and thus one can use either one of the two notions freely.

\section{The discrepancy of subsets of $\R^n$} \label{sec:disc}
We begin this section with a technical lemma which is at the
heart of the proof of Theorem \ref{theoA}. The lemma allows one to find a good
choice of signs on roughly half of the coordinates, while the error
incurred by the choice of coordinates and signs can be controlled
using the geometric structure of $T$.

A preliminary result we need has to do with Bernoulli processes---the
well-known H\"{o}ffding inequality (see, e.g., \cite{LT,VW}).
\begin{Theorem} \label{thm:hoeffding}
Let $(\eps_i)_{i=1}^n$ be independent, symmetric, $\{-1,1\}$-valued
random variables. Then, for every $a \in\R^n$ and every $t>0$,
\[
\Prr \Biggl(\sum_{i=1}^n \eps_i a_i \geq t|a|\Biggr) \leq\exp(-t^2/2).
\]
In particular,
\[
\Prr \Biggl(\Biggl| \sum_{i=1}^n \eps_i a_i \Biggr| \geq t|a| \Biggr)
\leq2\exp(-t^2/2).
\]
\end{Theorem}

Let us formulate the main lemma.
\begin{Lemma} \label{lemma:gen-disc}
Let
\[
\Phi(t)=
\cases{
\log(e/t), &\quad if $0<t \leq1$,\cr
t\exp(-t+1), &\quad if $t>1$.}
\]
There exist absolute constants $\ka_1$ and $\ka_2$
for which the following holds.
Assume that $(\lambda_s)_{s=1}^\infty$ is an increasing positive sequence
tending to infinity,
$(Q_s)_{s=1}^\infty$ is a positive sequence and $n$ is an integer such that
\[
\ka_1 \sum_{s=1}^\infty\lambda_s \Phi((\ka_2 Q_s)^2) \leq\frac{n}{100}.
\]
Let $T \subset\R^n$ for which $0 \in T$, set $(T_s)_{s=1}^\infty$
to be a sequence of subsets of $T$ and let $T_0=\{0\}$. Consider
maps $\pi_s\dvtx T \to T_s$ that satisfy that:
\begin{enumerate}
\item[(a)] for every $s \geq1$, $|\{\pi_s(t)-\pi_{s-1}(t) \dvtx t
\in T\}| \leq
\lambda_s$,
\item[(b)] for every $t \in T$, $\lim_{s \to\infty} \pi_s(t)=t$.
\end{enumerate}

Then, there exists $(\eta_i)_{i=1}^n \in\{-1,0,1\}^n$ such that
$n/4 \leq|\{i \dvtx \eta_i = 0\}| \leq3n/4$, and for every $t \in T$,
\[
\Biggl|\sum_{i=1}^n \eta_i t_i \Biggr| \leq\sum_{s=1}^\infty
Q_s|\pi_s(t)-\pi_{s-1}(t)|.
\]
\end{Lemma}

The proof is a combination of a chaining argument
and the entropy method, which is frequently used in
Discrepancy Theory (see, e.g., \cite{AS,Mat,Spen}).
In the chaining mechanism, one takes the sets $T_s$ to be finer and
finer approximations of the set $T$ and $\pi_s(t)$ is a nearest
element to $t$ in $T_s$ with respect to the underlying metric
(which is, in our case, the $\ell_2^n$ metric).

Recall that the entropy of a discrete random variable $X$ taking
values in a countable set $\Omega$ is
\[
H(X)=-\sum_{\omega\in\Omega} p_\omega\log_2 p_\omega,
\]
where $p_\omega= \Prr(X = \omega)$. The entropy function $H(X)$
indicates how close $X$ is to being equally distributed; the more
equally distributed $X$ is, the larger $H(X)$ is.

The three facts we will need regarding the entropy are well known
and we omit their proofs. First, if $H(X) \leq K$ then there is a
value of $X$ that is attained with probability at least $2^{-K}$.
Second, if $X$ attains at most $k$ values then $H(X) \leq\log_2 k$,
and finally, if $X=(X_1,\ldots,X_m)$ is a random vector then $H(X) \leq
\sum_{i=1}^m H(X_i)$.

In the entropic argument we will use, each ``link'' in each chain in
$T$ is assigned a random variable $X_\alpha\dvtx\{-1,1\}^n \to\R$ that
depends on the link and on the chain. The idea is to show that with
probability at least $2^{-\eta n}$, for every $\alpha$, each random
variable $X_\alpha$ falls in an interval $I_\alpha$ whose length is at
most $Q_\alpha\|X_\alpha\|_{L_2}$. One would like to make these
scaling factors $Q_\alpha$ as small as possible while still ensuring
that conditions $1$ and $2$ hold, since those conditions imply
that the intersection of the level sets of all the random variables
$X_\alpha$ has the desired measure.

More details on the way entropic arguments have been used in the
context of Discrepancy Theory may be found in \cite{AS,Matbook}.

Before presenting the proof, one should mention that a chaining
argument was implicit in Matou\v{s}ek's result on the discrepancy of
a subset of $\{0,1\}^n$ with a bounded VC dimension
\cite{Mat,Matbook}.

The first step in the proof of Lemma \ref{lemma:gen-disc} is the following
entropy estimate. We denote by $[x]$ the integer value of $x$.
\begin{Lemma} \label{lemma:entropic-estimate}
There exists an absolute constant $c$ for which the following holds.
Let $a \in\R^n$, set $Z_a = \sum_{i=1}^n \eps_i a_i$ and put
\[
W_a=\sgn(Z_a)[ | Z_a | ].
\]
Then
\[
-\sum_{i=-\infty}^\infty \Prr(W_a=i) \log \Prr(W_a=i) \leq c\Phi(1/2|a|^2).
\]
\end{Lemma}
\begin{pf}
By H\"{o}ffding's inequality
(Theorem \ref{thm:hoeffding}), for every $j \in\mathbb{Z} \setminus
\{0\}$,
\[
p_j = \Prr ( W_a = j )
\leq \Prr ( Z_a \geq|j| ) < \exp(-j^2/2|a|^2).
\]
Also,
\[
p_0=\Prr (W_{a} = 0 )=\Prr (-1 < Z_{a} < 1
),
\]
implying that
\[
1-p_0 =\Prr \Biggl( \Biggl|\sum_{i=1}^n \eps_i a_i \Biggr| \geq1 \Biggr) \leq
2\exp(-1/2|a|^2).
\]
Consider $j \in\Z$ for which $|j| \geq\sqrt{2}|a|$. Since
$f(x)=-x\log x$ is increasing in
$[0,1/e]$, it follows that for such values of $j$,
\[
-p_j \log p_j \leq\frac{j^2}{2|a|^2} \exp(-j^2/2|a|^2).
\]
Fix an integer $k$ which satisfies that $k \geq\sqrt{2}|a|$ and which
will be named later, and observe
that if we set $S=\sum_{1 \leq|j| \leq k} p_j$ then
\[
-\sum_{1 \leq|j| \leq k} p_j \log p_j \leq-\sum_{1 \leq|j| \leq k}
\frac{S}{2k} \log(S/2k) = S\log(2k/S).
\]
Clearly, $S \leq1-p_0 \leq2\exp(-1/2|a|^2)$, and thus, if $\exp
(-1/2|a|^2) \leq1/e$ (i.e., if $\sqrt{2}|a| \leq1$), then
\[
S\log(2k/S)= \log k + 2(S/2)\log(2/S) \leq\log k + \frac{1}{2|a|^2}
\exp(-1/2|a|^2).
\]
Otherwise, $S\log(2k/S) \leq\log k + \frac{1}{e}\log(2e) \leq
1+\log k $, implying that
\[
-\sum_{1 \leq|j| \leq k} p_j \log p_j \leq\log k +
\cases{\displaystyle
\frac{1}{2|a|^2} \exp(-1/2|a|^2), &\quad if $\sqrt{2}|a| \leq1$,\vspace*{2pt}\cr
1, &\quad otherwise.}
\]
Moreover,
\begin{eqnarray*}
-\sum_{|j| \geq k+1} p_j \log p_j &\leq& 2 \sum_{j \geq k+1} \frac
{j^2}{2|a|^2} \exp(-j^2/2|a|^2)
\\
&\leq& 2\int_k^\infty\frac{x^2}{2|a|^2}\exp(-x^2/2|a|^2) \\
&\leq&
(k+2|a|)\exp(-k^2/2|a|^2).
\end{eqnarray*}

Therefore,
\begin{eqnarray*}
-\sum_{j=\infty}^\infty p_j \log p_j & = & - \sum_{|j| > k } p_j \log
p_j - p_0\log p_0 - \sum_{1 \leq|j| \leq k } p_j \log p_j
\\
& \leq & (k+2|a|)\exp(-k^2/2|a|^2) + 2\exp(-1/2|a|^2)
\\
&&{} + \log k +
\cases{\displaystyle
\frac{1}{2|a|^2} \exp(-1/2|a|^2), &\quad if $\sqrt{2}|a| \leq1$,\vspace*{2pt}\cr
1, &\quad otherwise.}
\end{eqnarray*}
Now, consider the following three cases. First, if $\sqrt{2}|a| \leq
1$, take $k=1$,
and thus
\[
-\sum_{j=\infty}^\infty p_j \log p_j \leq\frac{c_1}{|a|^2}\exp(-1/2|a|^2).
\]
If $1<\sqrt{2}|a|\leq e$ set $k$ to be a suitable absolute constant and
if $\sqrt{2}|a| > e$, put $k \sim|a|\log(\sqrt{2}|a|)$. Therefore,
in both these cases
\[
-\sum_{j=\infty}^\infty p_j \log p_j \leq{c_2}\log\bigl(\sqrt{2}|a|\bigr),
\]
and our claim follows.
\end{pf}
\begin{pf*}{Proof of Lemma \ref{lemma:gen-disc}}
Without loss of
generality, assume that $T$ is finite. Recall that $0 \in T$ and that
$T_0=\{0\}$, consider the sets $T_s$ and the maps $\pi_s\dvtx T \to T_s$,
let $\Delta_s(t)=\pi_s(t)-\pi_{s-1}(t)$ and put $\Delta_s =
\{\pi_{s}(t)-\pi_{s-1}(t) \dvtx t \in T\}$. Let
$(\lambda_s)_{s=1}^\infty$ and $(Q_s)_{s=1}^\infty$ be as in
the assumptions of the lemma and set $(\eps_i)_{i=1}^n$ to be
independent, symmetric, $\{-1,1\}$-valued random variables.

Consider the Bernoulli process $t \to Z_t = \sum_{i=1}^n \eps_i t_i$.
Since $Z_t$ is linear in $t$ and $\pi_0(t)=0$, then for every $t \in
T$,
\[
t=\sum_{s=1}^\infty\Delta_s(t) \quad\mbox{and}\quad  Z_t
=\sum_{s=1}^\infty Z_{\pi_s(t)} -Z_{\pi_{s-1}(t)} =
\sum_{s=1}^\infty Z_{\Delta_s(t)}.
\]
For every $s \geq1$ and $u \in\Delta_s$ define
\[
\tilde{W}_{u,s} =\frac{Z_u}{|u|Q_s},\qquad W_{u,s} = \sgn
(\tilde{W}_{u,s}) [|\tilde{W}_{u,s}|].
\]

Observe that $(W_{u,s})_{u \in\Delta_s, s=1,2,\ldots}$
is a vector that takes a finite number of values. Since the entropy is
subadditive then
\[
H\bigl((W_{u,s}) \dvtx u \in\Delta_s, s=1,2,\ldots\bigr)
\leq\sum_{s=1}^\infty\sum_{u \in\Delta_s} H(W_{u,s}) =(*).
\]
Suppose that one can find $(Q_s)_{s=1}^\infty$ for which $(*) \leq
n/100$. By the properties of the entropy, this implies that there
are numbers $\{\ell_{u,s} \in\Z\dvtx u \in\Delta_s, s=1,2,\ldots\}$ such
that
%
\begin{equation} \label{eq:A}
\Prr \bigl( (\eps_i)_{i=1}^n\dvtx\forall u \in\Delta_s,
s\geq1,W_{u,s} = \ell_{u,s} \bigr) \equiv \Prr(A) \geq2^{-n/100}.
\end{equation}
Since $|A| \geq2^{0.99n}$, there will be at least two vectors
$(\eps_i)_{i=1}^n$ and $(\eps_i^\prime)_{i=1}^n$ in $A$ that differ
on at most $3n/4$ coordinates and on at least $n/4$ of them. The
desired sequence will then be
$(\eta_i)_{i=1}^n=(\frac{\eps_i-\eps_i^\prime}{2})_{i=1}^n$. Indeed,
for $u \in\Delta_s$,
\[
\Biggl|\sum_{i=1}^n \eta_i u_i \Biggr| =\frac{1}{2}\Biggl|\sum_{i=1}^n
\eps_i u_i - \sum_{i=1}^n \eps_i^\prime u_i \Biggr| \leq Q_s|u|,
\]
implying that every $t \in T$ satisfies
\[
\Biggl|\sum_{i=1}^n \eta_i t_i\Biggr| = \Biggl|\sum_{s \geq1}
\sum_{i=1}^n \eta_i (\Delta_s(t))_i \Biggr| \leq
\sum_{s=1}^\infty Q_s |\Delta_s(t)|.
\]

Hence, to complete the proof, it remains to show that for a sequence
$(Q_s)_{s=1}^\infty$ that satisfies the assumptions of the lemma,
$(*) \leq n/100$. Applying Lemma \ref{lemma:entropic-estimate} for
$a=u/|u|Q_s$, and since $1/2|a|^2=Q_s^2/2$, it is evident that
$H(W_{u,s}) \lesssim\Phi(Q_s^2/2)$. Thus,
\[
\sum_{s=1}^\infty\sum_{u \in\Delta_s} H(W_{u,s}) \lesssim
\sum_{s=1}^\infty|\Delta_s| \sup_{u \in\Delta_s} H(W_{u,s})
\lesssim\sum_{s=1}^\infty\lambda_s \Phi(Q_s^2/2),
\]
proving our claim.
\end{pf*}

We will apply Lemma \ref{lemma:gen-disc} in two typical situations.
The first case will lead to a bound on the discrepancy of a set
using the $\gamma_{2,s}$ functionals of the set and of its
coordinate projections. The second will result is an entropy
integral type bound, presented in Section \ref{sec:entropy-integral},
which will then be used to re-prove Spencer's
result on the discrepancy of a finite set system
\cite{Spen,Matbook} and Matou\v{s}ek's VC theorem
\cite{Mat,Matbook}.

Corollary \ref{cor:choice-of-Q-for-gamma-2} below will play a
central part in the proof of Theorem \ref{theoA}. Since it follows from a simple
computation, we omit its proof.
\begin{Corollary} \label{cor:choice-of-Q-for-gamma-2}
There exist absolute constants $\ka_3$, $\ka_4$, $\ka_5$
for which the following holds. Let $T \subset\ell_2^n$, assume that
$0 \in T$,
set $s_n=\max\{s\dvtx 2^{2^{s+1}} \leq\ka_3 n\}$ and put $T_s$ to be a
collection of subsets of $T$ with $|T_s| \leq2^{2^s}$. Then, if
\[
Q_s = \ka_4
\cases{
\exp(-\ka_5n^{1/2}), &\quad if $s<s_n$,\cr
1, &\quad if $s=s_n$,\cr
2^{s/2}, &\quad if $s>s_n$,}
\]
there exists $(\eta_i)_{i=1}^n \in\{-1,0,1\}^n$ such that $n/4 \leq
|\{i \dvtx \eta_i=0\}| \leq3n/4$, and for every $t \in T$,
\[
\Biggl|\sum_{i=1}^n \eta_i t_i \Biggr| \leq\sum_{s=1}^\infty Q_s
|\pi_s(t)-\pi_{s-1}(t)|,
\]
where $\pi_s(t)$ is a nearest point to $t$ in $T_s$.
\end{Corollary}

\subsection{An entropy integral argument} \label{sec:entropy-integral}
In this section, we will prove an analog of Dudley's entropy integral
bound (see, e.g., \cite{LT,Tal-book}) in the context of discrepancy.
The entropy integral is often used to upper bound ${\sup_{t \in T}}
|{\sum_{i=1}^n \eps_i t_i}|$ for a typical $(\eps_i)_{i=1}^n$, but
here we will present a modified version that allows one to control
${\inf_{\eta} \sup_{t \in T}} |{\sum_{i=1}^n \eta_i t_i }|$, where the
infimum is taken with respect to all $\eta=(\eta_i)_{i=1}^n \in
\{-1,0,1\}^n$ for which roughly half the coordinates are nonzero.

Let $T \subset\ell_2^n$ and recall that for every $\eps>0$,
$D(\eps)=D(\eps,T,\ell_2^n)$ is the cardinality of a maximal
$\eps$-separated subset of $T$. Also, set
\[
u(\eps)=
\cases{
\sqrt{\displaystyle \log\biggl(\frac{eD(\eps)}{n}\biggr)}, &\quad if $D(\eps) \geq
n$,\vspace*{2pt}\cr
\displaystyle \exp\Biggl(-\sqrt{\frac{n}{D(\eps)}}+1\Biggr), &\quad if $D(\eps)<n$.}
\]
\begin{Theorem} \label{thm:entropy-integral-disc}
There exist an absolute constant $c$ for
which the following holds. If $T \subset\ell_2^n$ and $0 \in T$, then
there exist $(\eta_i)_{i=1}^n \in\{-1,0,1\}^n$,
such that $n/4 \leq|\{i \dvtx \eta_i = 0\}| \leq3n/4$
and for every $t \in T$,
%
\begin{equation}\label{eq:entropy-integral-disc}
\Biggl|\sum_{i=1}^n \eta_i t_i \Biggr| \leq c \int_0^{\operatorname{diam}(T)} u(\eps)\,d\eps.
\end{equation}
\end{Theorem}
\begin{Remark}
Recall that Dudley's entropy integral bound shows that
\[
\E\sup_{t \in T} \Biggl|\sum_{i=1}^n \eps_i t_i \Biggr|
\leq c_1 \int_0^{\operatorname{diam}(T)} \sqrt{\log D(\eps)} \,d\eps,
\]
for a suitable absolute constant $c_1$. Clearly, this entropy integral may
be considerably larger than the quantity we have in Theorem
\ref{thm:entropy-integral-disc}. It is also evident that if one
could iterate Theorem \ref{thm:entropy-integral-disc} for the set
$P_I T \subset\ell_2^{|I|}$, where $I=\{i\dvtx \eta_i=0\}$, and continue
in the same manner, then one would likely improve upon the bound
resulting from the standard entropy integral bound that holds for
a typical choice of signs, if indeed distances in $P_I T$ shrink
relative to distances in $T$.
\end{Remark}

The proof of Theorem \ref{thm:entropy-integral-disc} is based on
Lemma \ref{lemma:gen-disc}. It
requires two additional simple results. Since their proofs are
standard, we shall not present them here.
\begin{Lemma} \label{lemma:q-est-1}
There exist absolute constants $c_1$, $c_2$, $c_3$ and $c_4$ for which the
following holds.
Let $T \subset\ell_2^n$, set $\nu_n$ to be the largest integer $s$
satisfying $2^s \leq c_1 n$ and define
\[
\lambda_s =
\cases{
c_2 2^s, &\quad if $s \leq\nu_n$,\cr
c_3 n 2^{2^{s-\nu_n}-1}, &\quad if $s>\nu_n$.}
\]
Then conditions \textup{(a)} and \textup{(b)} of Lemma
\ref{lemma:gen-disc} hold if one selects
\[
Q_s = c_4
\cases{
\exp\bigl(-2 \cdot2^{(s-\nu_n)/2}\bigr), &\quad if $s \leq\nu_n$,\vspace*{1pt}\cr
2^{(s-\nu_n)/2}, &\quad if $s > \nu_n$.}
\]
\end{Lemma}
\begin{Lemma} \label{lemma:integral-vs-sum}
Let $g$ and $f$ be nonincreasing, nonnegative functions and let
$(\eps_s)_{s=0}^m$ be a decreasing sequence. If for every $s \geq
1$, $g(\eps_{s-1}) \geq f(\eps_s)$, and if there is $\alpha>0$ such
that for every $s \geq1$, $f(\eps_s)-f(\eps_{s-1}) \geq\alpha
f(\eps_s)$ then
\[
\int_{\eps_m}^{\eps_0} g(\eps)\,d\eps+ \eps_m f(\eps_m) \geq
\alpha
\sum_{s=1}^m f(\eps_s)\eps_{s-1}.
\]
\end{Lemma}
\begin{pf*}{Proof of Theorem \ref{thm:entropy-integral-disc}}
Let $(\lambda_s)_{s=1}^\infty$ and $(Q_s)_{s=1}^\infty$ be as in
Lemma \ref{lemma:q-est-1}. Without loss of generality assume that
$T$ is a finite set and define the sets $T_s$ iteratively, as
follows. Set $m$ to be the first integer such that $|T| \leq
\lambda_m$, let $T_s=T$ for $s \geq m$ and set $\eps_m=0$.
For $m-1$, let
$\eps_{m-1}=\inf\{\eps\dvtx D(\eps,T_m,\ell_2^n) \leq\lambda_{m-1}\}$
and put $T_{m-1}$ to be a maximal $\eps_{m-1}$-separated subset
of $T_m$ whose cardinality is at most $\lambda_{m-1}$.
Continue in this way to construct the sets $T_s$ for $s=m-1,\ldots,1$.
For every~$s$, let $\pi_s(t)$ be a nearest point to $\pi_{s+1}(t)$ in $T_s$.

Let $s \leq m$, and since the sets $T_s$ are nested, then
$|\{\pi_s(t)-\pi_{s-1}(t) \dvtx t \in T\}| \leq|T_s|$ and
$|\pi_s(t)-\pi_{s-1}(t)| \leq\eps_{s-1}$ for every $t\in T$. Therefore,
applying Lemmas \ref{lemma:gen-disc} and~\ref{lemma:q-est-1}, there is
a choice $(\eta_i)_{i=1}^n \in\{-1,0,1\}^n$ with $n/4 \leq|\{i \dvtx
\eta_i = 0\}| \leq3n/4$ such that for every $t \in T$,
%
\begin{equation} \label{eq:upper.sums}
\Biggl| \sum_{i=1}^n \eta_i t_i \Biggr| \leq c_1 \biggl(\sum_{s
\leq
\nu_n} \exp\bigl(-2 \cdot2^{(\nu_n-s)/2}\bigr)\eps_{s-1} + \sum_{s>\nu_n}^{m}
2^{(s-\nu_n)/2}\eps_{s-1} \biggr).
\end{equation}

It remains to bound the sums in (\ref{eq:upper.sums}) by the
appropriate integrals,
using Lem\-ma~\ref{lemma:integral-vs-sum}. First, for $s > \nu_n$ let
\[
f(\eps)= \sum_{s=\nu_n+1}^{m} 2^{(s-\nu_n)/2}
\IND_{(\eps_{s+1},\eps_s]},\qquad  g(\eps)= \sum_{s=\nu_n+1}^{m}
2^{(s-\nu_n)/2} \IND_{(\eps_{s},\eps_{s-1}]}.
\]
Clearly, in $[\eps_{m},\eps_{\nu_n}]=[0,\eps_{\nu_n}]$, $f$ and $g$
are nonincreasing and nonnegative, for every $\eps$ in that range
\[
f(\eps) \leq g(\eps) \leq\sqrt{2}f(\eps) \lesssim u(\eps),
\]
and the conditions of Lemma \ref{lemma:integral-vs-sum} hold. Since
$\eps_{m}=0$, then
\[
\sum_{s>\nu_n}^{m} 2^{(s-\nu_n)/2}\eps_{s-1} \lesssim
\int_0^{\eps_n} u(\eps).
\]
For the other term in (\ref{eq:upper.sums}), if $s \leq\nu_n$ then
$2^s \sim\lambda_s \leq D(\eps_s)$ and the sum is estimated in a
similar way.
\end{pf*}

\subsubsection{Spencer's theorem}
Let us show how Theorem \ref{thm:entropy-integral-disc} can be used
to prove a version of Spencer's celebrated result from \cite{Spen}
(see also \cite{AS,Matbook}).
\begin{Theorem} \label{thm:Spencer}
There exists an absolute constant $c$ such that if
$T \subset B_\infty^n$ is of cardinality $m \geq n$, then
\[
\operatorname{disc}(T) \leq c \sqrt{n \log\biggl(\frac{em}{n}\biggr)}.
\]
\end{Theorem}
\begin{pf}
Without loss of generality, assume that $0 \in T$. Using the notation
of Theorem
\ref{thm:entropy-integral-disc}, for every $\eps>0$, $u(\eps) \leq
\sqrt{\log(
em/n)}$, and since $T \subset B_\infty^n$ then $\operatorname{diam}(T) \leq
\sqrt{n}$. Hence,
there are $(\eta_i)_{i=1}^n \in\{-1,0,1\}^n$ for which $n/4 \leq
|\{i \dvtx \eta_i=0\}| \leq3n/4$ and for every $t \in T$,
\[
\Biggl|\sum_{i=1}^n \eta_i t_i \Biggr| \leq c_1
\int_0^{\sqrt{n}} \sqrt{\log(e m/{n})} \,d\eps
\leq c_1\sqrt{n \log
({em}/{n})}.
\]
Now the result follows by repeating this argument for $P_{I_1} T$,
where $I_1 = \{i \dvtx \eta_i = 0\}$, an so on.
\end{pf}

\subsubsection{\texorpdfstring{Matou\v{s}ek's VC theorem}{Matousek's VC theorem}}
A well-known measure of complexity for subsets of $\{0,1\}^n$ is
the $\VC$ dimension of the set (its real value counterpart will be used
in Section \ref{sec:upp-low}).
\begin{Definition}
Let $T \subset\{0,1\}^n$. We say that $\sigma\subset\{1,\ldots,n\}$
is shattered by $T$ if $P_\sigma T = \{0,1\}^\sigma$---that is, if
the coordinate projection $P_\sigma T=\{(t_i)_{i \in\sigma} \dvtx t \in
T\}$ is the entire combinatorial
cube on these coordinates. Define $\VC(T)$ to be the maximal cardinality
of a subset of $\{1,\ldots,n\}$ that is shattered by $T$.
\end{Definition}

In \cite{Mat}, Matou\v{s}ek proved that the discrepancy of a VC
class is polynomially better than could be expected from a random
choice of signs. He obtained the best possible estimate for the
discrepancy of VC-subsets of $\{0,1\}^n$ as a function of the
dimension $n$.
\begin{Theorem} \label{thm:mat}
For every integer $d$, there is a constant $c(d)$ for which the
following holds. If $T \subset\{0,1\}^n$ and $\VC(T) \leq d$, then
$\operatorname{disc}(T) \leq c(d)n^{1/2-1/2d}$.
\end{Theorem}

To prove Matou\v{s}ek's theorem, recall the following fundamental
property of a VC class, due to Haussler \cite{Hau}.
\begin{Lemma} \label{lemma:vc-cover}
If $T \subset\{0,1\}^n$ and $\VC(T) \leq d$, then for every
$I \subset\{1,\ldots,n\}$ and every $0 < \eps\leq|I|^{1/2}$,
\[
D(\eps, P_I T, \ell_2^{I} ) \leq c(d)
\biggl(\frac{|I|^{1/2}}{\eps}\biggr)^{2d},
\]
where $c(d)$
is a constant that depends only on $d$.
\end{Lemma}
\begin{pf*}{Proof of Theorem \ref{thm:mat}}
Again, we may assume that $0 \in T$ and view $T$ as a subset of $\R^n$.
Let $\eps_n=\inf\{ \eps\dvtx D(\eps) \leq n\}$. Therefore, $\eps_n
\leq
c_1(d) n^{1/2-1/2d}$. A~change of variables shows that
\begin{eqnarray*}
\int_0^{\eps_n} \sqrt{\log\bigl({eD(\eps)}/{n}\bigr)} \,d\eps
& \leq c_2(d) n^{1/2-1/2d},
\\
\int_{\eps_n}^{\operatorname{diam}(T)} \exp\bigl(- \sqrt{{n}/{D(\eps
)}}\bigr)
\,d\eps& \leq c_2(d) n^{1/2-1/2d}.
\end{eqnarray*}
Hence, there is a choice of $(\eta_i^0)_{i=1}^n \in
\{-1,0,1\}^n$ such that for every $t \in T$
\[
\Biggl|\sum_{i=1}^n \eta_i^0 t_i \Biggr| \leq c_3(d)n^{1/2-1/2d},
\]
and if we set $I_1=\{i \dvtx \eta_i^1=0\}$ then $|I_1| \leq3n/4$. Since
$\VC(P_{I_1} T) \leq d$ then repeating the same argument for the set
$P_{I_1} T$, there are $(\eta^1_i)_{i \in I_1} \in\{-1,0,1\}^{I_1}$
such that for every $t \in T$, $|{\sum_{i \in I_1} \eta_i^1 t_i}
| \leq c_3(d)|I_1|^{1/2-1/2d}$, and so on. Therefore, there is a
choice of signs $(\eps_i)_{i=1}^n$ such that for every $t \in T$,
\[
\Biggl|\sum_{i=1}^n \eps_i t_i \Biggr| \leq c_3(d) \sum_j
|I_j|^{1/2-1/2d} \leq c_4(d) n^{1/2-1/2d},
\]
where we have used the fact that for every $j$, $|I_j| \leq
3|I_{j-1}|/4$.
\end{pf*}

The proof of Theorem \ref{thm:mat} illustrates once again the main
property we used to bound the discrepancy of a subset of $\R^n$. It
is not enough for the set to be small in the sense of its metric
entropy; what is needed is additional control on the ``size'' of all
of the set's coordinate projections. One way of controlling those
coordinate projections is
by taking into account information about the position of vectors
in the set, since coordinate projections of vectors in a good
position shrink norms and mutual distances.

\section{A decomposition theorem for subgaussian processes} \label{sec:shrink}
It is clear from our estimate on the discrepancy of a set $T \subset
\R^n$ that it would be useful to control the $\ell_2^I$ distances
between points in $T$ for every $I \subset\{1,\ldots,n\}$---that is,
distances between coordinate projections of elements of $T$. One
would be able to obtain a good bound on $\operatorname{disc}(T)$ if $T$ is
not too rich and if for every $I \subset\{1,\ldots,n\}$ and every $x,y
\in T$, $\|x-y\|_{\ell_2^{I}}$ is significantly smaller than
$\|x-y\|_{\ell_2^n}$. Unfortunately, usually this is not true even
for a single vector $z=x-y$. Indeed, if $z$ is supported in $I$ then
$\|z\|_{\ell_2^I}$ does not ``shrink'' at all. On the other hand, if
the coordinates of $z$ are roughly equal, then the coordinate
projection onto any $I$ shrinks the norm of $z$ by a factor of
$(|I|/n)^{1/2}$.

It is well known that a strong shrinking phenomenon is exhibited by
vectors in a general position. In this section, we will show that
if a class of functions $F$ is $L$-subgaussian, then a shrinking
phenomenon happens for a typical set
\[
T=P_\sigma F = \{ (f(X_i))_{i=1}^k \dvtx f \in F \},
\]
uniformly for all coordinate projections of $T$.

\subsection{Shrinking for a single function}
As a starting point, let us describe the so-called ``standard
shrinking'' phenomenon for a single function $f$. Let $f$ be a
function for which $\|f\|_{\psi_2} \leq L \|f\|_{L_2}$. Then,
concentration implies that with high probability,
\[
\|f\|_{L_2^k} = \Biggl(\frac{1}{k} \sum_{i=1}^k f^2(X_i)
\Biggr)^{1/2} \sim\|f\|_{L_2}.
\]
However, as we mentioned above, the shrinking phenomenon one needs
here is more general---that for every subset $I \subset
\{1,\ldots,k\}$, the $L_2^I$ norm of $f$ is upper bounded (possibly up
to a logarithmic factor) by $\|f\|_{L_2}$ (which translates in the
$\ell_2$ normalization to the shrinking of the norm). The following
lemma shows that this stronger claim is true as well whenever
$f$ is $L$-subgaussian.
\begin{Lemma} \label{lemma:shrinking-single}
For every $0<\delta<1$ and $L>0$, there is a constant $c(\delta,L)$
for which the following holds. If $\|f\|_{\psi_2} \leq L
\|f\|_{L_2}$ then for every integer $k$, with probability at least
$1-\delta$, for every $I \subset\{1,\ldots,k\}$,
\[
\|f\|_{L_2^I} \leq c(\delta,L) \sqrt{\log{({ek}/{|I|})}}
\|f\|_{L_2}.
\]
\end{Lemma}
\begin{pf}
Fix $k$ and $I\subset\{1,\ldots,k\}$. Since
$\|f^2\|_{\psi_1}=\|f\|_{\psi_2}^2$, then by Bernstein's inequality,
for every $t>0$,
\[
\Prr \biggl( \biggl| \frac{1}{|I|} \sum_{i \in I} f^2(X_i) - \E f^2
\biggr| \geq t\|f\|_{\psi_2}^2 \biggr) \leq2 \exp(-c_0|I|
\min\{t^2,t\}).
\]
Let $m \leq c_1 k$ and recall that there are at most $(ek/m)^m$
subsets of $\{1,\ldots,k\}$ of cardinality $m$. Hence, it suffices to
take $t=\beta(\delta) \log(ek/m) \geq1$ and obtain that with
probability of at least $1-2\exp(-c_0\beta m \log(ek/m))$, for every
subset $I$ of $\{1,\ldots,k\}$ of cardinality $m$,
%
\begin{equation} \label{eq:single-function}
\|f\|_{L_2^I}=\biggl(\frac{1}{m}\sum_{i \in I} f^2(X_i) \biggr)^{1/2}
\leq c_2(\delta)L\sqrt{\log{({ek}/{m})}} \|f\|_{L_2}.
\end{equation}
Therefore, summing the probabilities with respect to $m$, it follows
that for the correct choice of $\beta$, with probability at least
$1-\delta$, (\ref{eq:single-function}) is true for all subsets of
$\{1,\ldots,k\}$ of cardinality at most $c_1k$. The claim now easily
follows.
\end{pf}

\subsection{Shrinking for a class of functions}
When one attempts to generalize this simple shrinking argument to a
class of functions, one faces a problem: the probabilistic estimate
obtained in the proof of Lemma \ref{lemma:shrinking-single} does not
allow one to control many functions simultaneously. Thus, a naive
extension of that result is simply too weak to lead to a
function class analog of the shrinking phenomenon.

To formulate the shrinking phenomenon for an $L$-subgaussian class
of functions, let us recall some notation. For any two sets $A$ and $B$
in a vector
space, $A+B = \{a+b\dvtx a \in A, b \in B\}$, and for a class of
functions $F$, a random sample $\sigma=(X_1,\ldots,X_k)$ and $I \subset\{
1,\ldots,k\}$,
\[
P_\sigma F = \{(f(X_i))_{i=1}^k \dvtx f \in F \},\qquad P_I^\sigma F =
\{(f(X_i))_{i \in I} \dvtx f \in F \}.
\]
For every integer $m$, let $W_m =\{x \in\R^m \dvtx x_j^* \leq
1/\sqrt{j}, j=1,\ldots,m \}$, where $(x_j^*)_{j=1}^m$ is a
monotone nonincreasing rearrangement of $(|x_j|)_{j=1}^m$. Thus,
$W_m$ is the unit ball of the weak $\ell_2$ space
$\ell_{2,\infty}^m$. Denote by $V_m$ the collection of all subsets
of $\{1,\ldots,k\}$ of cardinality at most $m$ and set $\tau_m$ to be
the smallest integer $s$ such that $2^{2^{s}} \geq\exp(m\log(ek/m))
\geq|V_m|$.
\begin{Theorem} \label{thm:decomposition-of-F}
For every $0<\delta<1$ and $L>0$ there exist constants $c_1$, $c_2$,
$c_3$ and $k_0$ depending only on $L$ and $\delta$ for which the
following holds. Let $F$ be an $L$-subgaussian class of functions and
assume that for each $f \in F$, $\E f =\alpha$ for some $\alpha\in
\R$. Then, for\vspace*{1pt} every integer $k$ and every $m \leq k$, there are
sets $F_1^{m}$ and $F_2^{m} \subset F$ with the following
properties. First, $F \subset F_1^{m} + F_2^{m}$; second, with
$\mu^k$-probability of at least $1-\delta$, if
$\sigma=(X_1,\ldots,X_k)$ then:
\begin{enumerate}
\item For every integer $m \leq k$ and every $I \subset\{1,\ldots,k\}$
of cardinality $m$,
\[
P_I^\sigma F_1^{m} \subset c_1 \gamma_{2,\tau_m}(F,L_2) W_{m}.
\]
\item For every $f,h \in F_2^{m}$ and every $I \subset
\{1,\ldots,k\}$ of cardinality $m$,
\[
\|f-h\|_{L_2^I} \leq c_2 \sqrt{\log(ek/m)} \|f-h\|_{L_2}.
\]
\item If $k \geq k_0$ then for every $m \leq c_3 k$ and
every $f,h \in F_2^m$,
\[
\|f-h\|_{L_2} \leq\sqrt{2} \|f-h\|_{L_2^\sigma}.
\]
\end{enumerate}
\end{Theorem}

The way Theorem \ref{thm:decomposition-of-F} should be understood is
as follows. Consider a typical $\sigma=(X_1,\ldots,X_k)$ and let $T
=P_\sigma F \subset\ell_2^k$. Then, for every $I \subset
\{1,\ldots,k\}$ the further coordinate projection satisfies $P_I T
\subset P_I T_1 + P_I T_2$ where $T_1,T_2 \subset\ell_2^k$ depend
only on the cardinality of $I$ and not on $I$ itself, and $T_2
\subset T$. The set $P_I T_1$ captures the ``peaky'' part of $P_I T$
and is contained in a relatively small set: a ball in
$\ell_{2,\infty}$ whose radius depends on the ``complexity'' of the
class $F$. The set $T_2$ consists of vectors that satisfy the
desired shrinking property. Indeed, for every $(f(X_i))_{i=1}^k,
(h(X_i))_{i=1}^k \in T_2$ and every $I \subset\{1,\ldots,k\}$ of
cardinality $m$ one has
\begin{eqnarray*}
\biggl(\sum_{i \in I} \bigl(f(X_i)-h(X_i)\bigr)^2 \biggr)^{1/2} & \leq & c_1
\sqrt{m \log({ek}/{m})} \|f-h\|_{L_2}
\\
& \leq & c_2 \sqrt{\frac{m}{k} \log({ek}/{m})}
\Biggl(\sum_{i=1}^k \bigl(f(X_i)-h(X_i)\bigr)^2 \Biggr)^{1/2},
\end{eqnarray*}
where the last inequality holds if $m \leq c_3k$.
\begin{pf*}{Proof of Theorem \ref{thm:decomposition-of-F}}
Fix an integer $k$. For every integer $m \leq k$, let
$(H_{s,m})_{s = \tau_m}^\infty$ be an almost optimal admissible
sequence of $F$ with respect to $\gamma_{2,\tau_m}(F,\psi_2)$, and
set $\pi_s^m$ to be the metric projection onto $H_{s,m}$ with
respect to the $\psi_2$ norm. For every such $m$ we will construct
two sets of functions, $F_1^{m}$ and $F_2^{m}$ such that $F \subset
F_1^{m} + F_2^{m}$ as follows: let $F_1^{m} = \{ f -
\pi^m_{\tau_m}(f) \dvtx f \in F\}$ and set $F_2^m=\{\pi^m_{\tau_m}(f) \dvtx f
\in F\}$ [and from here on we will omit the superscript $m$ and
write $\pi_s(f)$ instead of $\pi_s^m(f)$]. Note that this choice of
decomposition depends only on $m$ and does not depend on $k$.

For every $I \in V_m$ set $Z_f^I =\sum_{i \in I} (f(X_i)-\E
f)$ and observe that
\[
Z_f^I -Z_{\pi_{\tau_m}(f)}^I = \sum_{s >\tau_m} Z_{\pi_{s}(f)}^I
-Z_{\pi_{s-1}(f)}^I = \sum_{s >\tau_m} \sum_{i \in I}
\bigl(\pi_s(f)-\pi_{s-1}(f)\bigr)(X_i),
\]
since the expectation of all the functions in $F$ is the same. Thus,
for every $f \in F$, $\pi_s(f)-\pi_{s-1}(f)$ has mean zero, and for
every $s > \tau_m$ and every $t \geq1$,
\begin{eqnarray*}
&&
\Prr \biggl(\biggl|\sum_{i \in I} \bigl(\pi_s(f)
-\pi_{s-1}(f)\bigr)(X_i)\biggr| \geq t\|\pi_s(f)
-\pi_{s-1}(f)\|_{\psi_2} \sqrt{|I|} \biggr)
\\
&&\qquad\leq 2\exp(-c_0t^2).
\end{eqnarray*}
Let $t=u2^{s/2}$ for $u \geq c_1$, where $c_1$ is a constant to
be named later. Because of our choice of $s$, $|V_m| \leq
2^{2^s}$ and $|H_{s,m}| \cdot|H_{s-1,m}| \leq2^{2^{s+1}}$, and thus
\begin{eqnarray*}
&& \Prr \bigl( \exists f \in F, I \in V_m\dvtx \bigl|Z_{\pi_{s}(f)}^I - Z_{\pi_{s-1}(f)}^I\bigr|
\geq u2^{s/2} \|\pi_s(f) -\pi_{s-1}(f)\|_{\psi_2} \sqrt{|I|}
\bigr) \\
&&\qquad \leq2^{2^{s+1}} |V_m| \cdot2\exp(-c_0u^22^s) \leq
\exp(-c_2u^22^s).
\end{eqnarray*}
Hence, summing over $s>\tau_m$, it follows that with probability at
least
\[
1-\sum_{s
>\tau_m}\exp(-c_2u^22^{s}) \geq1-\exp(-c_3u^22^{\tau_m}),
\]
for every $f \in F$ and every $I \in V_m$
\[
\biggl| \sum_{i \in I} \bigl(f - \pi_{\tau_m}(f)\bigr)(X_i)
\biggr|
\leq u\sqrt{|I|} \sum_{s > \tau_m} 2^{s/2}\|\pi_s(f)
-\pi_{s-1}(f)\|_{\psi_2}.
\]
Summing the probabilities for all possible integers $1 \leq m \leq
k$ and noting that for every $1 \leq m \leq k$, $2^{\tau_m} \gtrsim
m\log(ek/m)$, it is evident that for $u \geq c_1$ there is a set
$\mathcal{A} \subset\Omega^k$ with probability at least
$1-\exp(-c_4u^2)$ for which the following holds. For every
$(X_i)_{i=1}^k \in\mathcal{A}$, every $1 \leq m \leq k$, every $h \in
F_1^m$ and every $I \in V_m$
\[
\biggl| \sum_{i \in I} h(X_i) \biggr| \leq2u \sqrt{|I|}
\gamma_{2,\tau_m}(F,\psi_2),
\]
where we have used the fact that $(H_{s,m})_{s = \tau_m}^\infty$ is
an almost optimal admissible sequence with respect to
$\gamma_{2,\tau_m}(F,\psi_2)$.

Fix $(X_1,\ldots,X_k) \in\mathcal{A}$, $1 \leq m \leq k$, $I \subset
\{1,\ldots,k\}$ of cardinality $m$ and $h \in F_1^m$. Consider the sets
$I_+=\{i \in I\dvtx h(X_i) \geq0\}$ and $I_{-}=\{i \in I\dvtx h(X_i) < 0\}$
and note that both are in $V_m$. Since $x^{1/2}$ is increasing, then
on the set $\mathcal{A}$
%
\begin{equation} \label{eq:f-1}
\sum_{i \in I} | h(X_i) | \leq4u \sqrt{|I|}
\gamma_{2,\tau_m}(F,\psi_2).
\end{equation}
In particular, if $(h_i^*)_{i=1}^k$ is a nonincreasing
rearrangement of $(|h(X_i)|)_{i=1}^k$ then by (\ref{eq:f-1}) applied
to the set $I_j$ consisting of the $j \leq m$ largest elements of
$(|h(X_i)|)_{i=1}^k$,
\[
h_{j}^* \leq\frac{1}{j} \sum_{i=1}^j h_i^* \leq\frac{1}{j} \cdot
4uj^{1/2}\gamma_{2,\tau_m}(F,\psi_2) \leq4Lu
\gamma_{2,\tau_m}(F,L_2)/\sqrt{j};
\]
thus, $P_I^\sigma F_1^m \subset4Lu \gamma_{2,\tau_m}(F,L_2) W_m$.

Turning our attention to the sets $F_2^m$, we will show that with
high probability, for every $1 \leq m \leq k$ and every $I \subset
\{1,\ldots,k\}$ of cardinality $m$, the coordinate projection $P_I
\dvtx(F_2^m,L_2) \to(F_2^m,L_2^I)$ has a well behaved Lipschitz
constant. To that end, fix $1 \leq m \leq k$, set $G_m=\{|f_1-f_2| \dvtx f_i \in F_2^m\}$ and recall that for every function $g$,
$\|g^2\|_{\psi_1} = \|g\|^2_{\psi_2}$. Hence, by Bernstein's
inequality, for every $g \in G_m$ and every $t>1$,
\[
\Prr \Biggl( \Biggl| \frac{1}{m}\sum_{i=1}^m g^2(X_i) -\E g^2 \Biggr|
\geq t \|g\|_{\psi_2}^2 \Biggr) \leq2\exp(-c_0 m \min(t^2,t)).
\]
Let $E_m$ be the collection of subsets of $\{1,\ldots,k\}$ of
cardinality $m$. Since $|G_m| \leq|F_2^m|^2 \leq2^{2^{\tau_m+1}}$
and $|E_m| \leq|V_m| \leq2^{2^{\tau_m}}$, then by taking $u \geq
c_5$ and $t=u\log(ek/m) \geq1$,
\begin{eqnarray*}
&&\Prr \biggl( \exists g \in G_m, I \in E_m \dvtx \biggl|
\frac{1}{m}\sum_{i \in I} g^2(X_i) -\E g^2 \biggr| \geq
\|g\|_{\psi_2}^2 u\log(ek/m) \biggr)
\\
&&\qquad\leq 2^{2^{\tau_m+2}} \exp\bigl(-c_0u m\log(ek/m)\bigr) \leq\exp\bigl(-c_6u
m\log(ek/m)\bigr).
\end{eqnarray*}
Summing over all possible $1 \leq m \leq k$, there is a subset
$\mathcal{B} \subset\Omega^k$ of probability at least
$1-\exp(-c_{7}u)$
on which the following holds. For every $1 \leq m \leq k$, every
$f_1,f_2 \in F_2^m$ and every $I \in E_m$,
%
\begin{eqnarray} \label{eq:f-2}
\|f_1-f_2\|_{L_2^I}^2 & \leq & \|f_1-f_2\|_{L_2}^2 + u
\log({ek}/{m}) \|f_1-f_2\|_{\psi_2}^2
\nonumber\\[-8pt]\\[-8pt]
& \leq & 2L^2 u \log({ek}/{m}) \cdot\|f_1-f_2\|_{L_2}^2.\nonumber
\end{eqnarray}

Thus, fix a ``legal'' choice of $u$ for which $\Prr(\mathcal{A} \cap
\mathcal{B}) \geq1 -\delta/2$. Since both (\ref{eq:f-1}) and (\ref{eq:f-2})
hold on that event, the proof of the first and second claims is
evident.

For the third part, fix $t <1/2 $ to be named later. Again, by
Bernstein's inequality and since $F$ is $L$-subgaussian, then with
probability at least $1-2|F_2^m|^2\exp(-c_0kt^2)$, for every
$f_1,f_2 \in F_2^m$
\[
\|f_1-f_2\|_{L_2}^2 \leq\|f_1-f_2\|_{L_2^k}^2 +
tL^2\|f_1-f_2\|_{L_2}^2.
\]
Thus, taking $t=1/(2L^2)$, for $k \geq k_0(\delta,L)$ and $m \leq
c_{8}(L)k$, it is evident that with probability at least
$1-2\exp(-c_{9}(L)k) \geq1-\delta/2$, for every $f_1,f_2 \in
F_2^m$,
\[
\|f_1-f_2\|_{L_2}^2 \leq2\|f_1-f_2\|_{L_2^k}^2,
\]
as claimed.
\end{pf*}

\subsection{Shrinking properties of the $\gamma_{2,s}$ functionals}

The first corollary of Theorem \ref{thm:decomposition-of-F} we shall
present here is a shrinking property of $\gamma_{2,s}(F,L_2^I)$.
\begin{Theorem} \label{thm:herd}
For every $0<\delta<1$ there exists a constant $c(\delta) \sim
\log(2/\delta)$ for which the following holds. Let $F$ be an
$L$-subgaussian class of functions on a probability space
$(\Omega,\mu)$ and assume that for every $f \in F$, $\E f = \alpha$
for some $\alpha\in\R$. Then, with probability at least
$1-\delta$, for every $I \subset\{1,\ldots,k\}$ and every integer $s$
that satisfies $2^s \leq|I| \log({ek/|I|})$,
\[
\gamma_{2,s+1}(F,L_2^I) \leq c(\delta)L \gamma_{2,s}(F,L_2)
\sqrt{\log({ek}/{|I|})}.
\]
\end{Theorem}

Before proving Theorem \ref{thm:herd}, recall the following
well-known result on the expectation of a monotone rearrangement of
independent standard Gaussian variables (see, e.g.,
\cite{GLSW,GLMP}).
\begin{Lemma} \label{lemma:gauss-arr}
Let $(g_i)_{i=1}^n$ be independent standard Gaussian variables and
denote by $(g_i^*)_{i=1}^n$ the nonincreasing rearrangement of
$(|g_i|)_{i=1}^n$. Then,
\[
\E g_i^* \sim
\cases{
\sqrt{\log({2n}/{i})}, &\quad if $i \leq n/2$, \vspace*{2pt}\cr
\displaystyle  1-\frac{i}{n+1}, &\quad if $i > n/2$.}
\]
Moreover,
\[
\Biggl(\E\sum_{i=1}^m (g_i^*)^2 \Biggr)^{1/2} \sim\sqrt{m \log
({en}/{m})}.
\]
\end{Lemma}
\begin{pf*}{Proof of Theorem \ref{thm:herd}}
Fix $0<\delta<1$ and let the sets $\mathcal{A},\mathcal{B} \in
\Omega^k$ be as in the proof of Theorem \ref{thm:decomposition-of-F}.
Take any $(X_1,\ldots,X_k) \in\mathcal{A} \cap\mathcal{B}$, let $I
\subset\{1,\ldots,k\}$ and set $m=|I|$. Since $P_I F \subset P_I F_1^m +
P_I F_2^m$, then by the sub-additivity of $\gamma_{2,s}$, it is evident
that for every integer $s$,
\[
\gamma_{2,s+1}(F,L_2^I) \leq
\gamma_{2,s}(F_1^m,L_2^I)+\gamma_{2,s}(F_2^m,L_2^I).
\]
By (\ref{eq:f-2}), the mapping $P_I \dvtx(F_2^m,L_2) \to(F_2^m,L_2^I)$
is a Lipschitz function with a constant
$c(\delta)L(\log(ek/m))^{1/2}$. Therefore, recalling that $F_2^m
\subset F$,
%
\begin{eqnarray} \label{eq:F_2}
\gamma_{2,s}(F_2^m,L_2^I) & \leq & c_1\sqrt{\log(ek/m)}
\gamma_{2,s}(F_2^m,L_2) \nonumber\\[-8pt]\\[-8pt]
& \leq & c_1\sqrt{\log(ek/m)} \gamma_{2,s}(F,L_2),\nonumber
\end{eqnarray}
where $c_1=c_1(L,\delta) \sim L \log(2/\delta)$. To conclude the proof,
observe that by Theorem \ref{thm:decomposition-of-F}, $P_I F_1^m
\subset B_m W_m$, where\vspace*{1pt}
$B_m=c_2(L,\delta)\gamma_{2,\tau_m}(F,L_2)$ and $W_m =\{x \in\R^m \dvtx
x_j^* \leq1/\sqrt{j},j=1,\ldots,m\}$.

Since the $\gamma_{2,s}$ functionals are monotone with respect to
inclusion and are decreasing in $s$, and since $\|x\|_{L_2^I} =
|I|^{-1/2}|x|$ for every $x \in\R^{m}$ then
\[
\gamma_{2,s}(F_1^m , L_2^I) \leq\gamma_2(F_1^m, L_2^I) \leq
B_m\frac{\gamma_2(W_m,| \cdot|)}{\sqrt{m}}.
\]
Applying the Majorizing Measures theorem and Lemma
\ref{lemma:gauss-arr}
\[
\gamma_2(W, | \cdot|) \leq c_{3} \E\sup_{w \in W} \sum_{i=1}^m g_i
w_i = c_{3} \E\sum_{i=1}^m \frac{g_i^*}{\sqrt{i}} \leq c_{4}
\sqrt{m}.
\]
Hence, for every $s$, $\gamma_{2,s}(F_1^m,L_2^I) \leq c_{4}B_m$,
implying that for every $s \leq\tau_m$, $\gamma_{2,s}(F_1^m,L_2^I)
\leq c_5(L,\delta) \gamma_{2,s}(F,L_2)$. Combining this with
(\ref{eq:F_2}), it follows that for every $I \subset\{1,\ldots,k\}$,
\[
\gamma_{2,s+1}(F,L_2^I) \leq c_6(L,\delta)
\sqrt{\log({ek}/{|I|})}\gamma_{2,s}(F,L_2),
\]
as claimed.
\end{pf*}
\begin{Remark} The proof of Theorem \ref{thm:herd}
yields a stronger result than the one formulated. It shows that with
probability $1-\delta$, for every $I \subset\{1,\ldots,k\}$ and every
$s \geq0$,
\[
\gamma_{2,s+1}(F,L_2^I) \leq
c(L,\delta)\bigl(\gamma_{2,\tau_{|I|}}(F,L_2)+\sqrt{\log(ek/|I|)}
\gamma_{2,s}(F,L_2)\bigr).
\]
Observe that in some sense, the range $s \leq\tau_{|I|}$ [i.e.,
$2^s \lesssim|I| \log(ek/|I|)$] is the interesting range of $s$, since
\[
\gamma_{2,s}(P_I^\sigma F , | \cdot|) \leq\operatorname{diam}(P_I^\sigma F
, | \cdot|) \gamma_{2,s} \bigl(B_2^{|I|},| \cdot|\bigr),
\]
which decreases exponentially in $s$ for $2^s \geq c_1|I|$.
\end{Remark}

Another outcome of Theorem \ref{thm:decomposition-of-F} was
formulated as Corollary \ref{corolB} in the \hyperref[intro]{Introduction}.
\begin{Corollary} \label{cor:Rn}
For every $0<\delta<1$ and $L>0$, there exist a constant
$c(\delta,L)$ such that the following holds. Let $\mu$ be an
isotropic, $L$-subgaussian measure on $\R^n$, set $(X_i)_{i=1}^k$ to
be independent, distributed according to $\mu$ and consider the
random operator $\Gamma=\sum_{i=1}^k \langle X_i,\cdot \rangle e_i$.
If $T
\subset\R^n$ and $V =k^{-1/2}\Gamma T$, then with
$\mu^k$-probability at least $1-\delta$, for every $I \subset
\{1,\ldots,k\}$,
\[
\E\sup_{v \in V} \biggl| \sum_{i \in I} g_i v_i \biggr| \leq
c(L,\delta) \sqrt{\frac{|I|}{k}\log({ek}/{|I|})} \E
\sup_{t \in T} \Biggl| \sum_{i=1}^n g_i t_i \Biggr|,
\]
where the expectation on both sides is with respect to the Gaussian
variables.
\end{Corollary}

The proof of Corollary \ref{cor:Rn} follows from Theorem
\ref{thm:decomposition-of-F} and the Majorizing Measures theorem.
\begin{pf*}{Proof of Corollary \ref{cor:Rn}}
Since $\mu$ is an $L$-subgaussian measure on $\R^n$, each $t
\in\R^n$ corresponds to a function $f_t(x) = \langle t,x \rangle$,
$f_t \dvtx \R^n \to\R$, for which $\|f_t\|_{\psi_2} \leq L|t|$. Let $F=\{f_t\dvtx
t \in T\}$, set $\Omega= \R^n$ and put $\sigma=(X_1,\ldots,\break X_k) \in
\Omega^k$ for which the assertion of Theorem
\ref{thm:decomposition-of-F} holds.

Fix $I \subset\{1,\ldots,k\}$ of cardinality $m$. Since $F$ is a class
of linear functionals, the decomposition of $F$ given in Theorem
\ref{thm:decomposition-of-F} actually implies a decomposition of $T$
which we denote by $T_1^m$ and $T_2^m$. Thus, for every $t \in T$,
$t=t^1+t^2$, where $t^i \in T_i^m$ for $i=1,2$. Since $ P_\sigma F
=\{ (f_t(X_i))_{i=1}^k \dvtx t \in T\} = \Gamma T$ then
\[
\E\sup_{v \in V} \biggl|\sum_{i \in I} g_i v_i \biggr| =
\frac{1}{\sqrt{k}} \E\sup_{t \in T} \biggl|\sum_{i \in I} g_i
\langle t_i^1+t_i^2,X_i \rangle \biggr|.
\]
Clearly, for any $u,v\in T$,
\[
\|f_u-f_v\|_{L_2^I} = m^{-1/2} \|(\langle X_i,u-v \rangle)\|_{\ell
_2^I}
\quad\mbox{and}\quad \|f_u-f_v\|_{L_2}=\|u-v\|_{\ell_2^n},
\]
and by the shrinking property of $T_2^m$, for every $u, v \in
T_2^m$,
\[
\|(\langle u,X_i \rangle)_{i=1}^n-(\langle v,X_i \rangle)_{i=1}^n\|
_{\ell_2^I} \leq
c_1(L,\delta) \bigl(m\log(ek/m)\bigr)^{1/2} \|u-v\|_{\ell_2^n}.
\]
Therefore, by Slepian's lemma (see, e.g., \cite{LT}) and since
$T_2^m \subset T$,
\begin{eqnarray*}
\frac{1}{\sqrt{k}} \E_g \sup_{t \in T} \biggl|\sum_{i \in I} g_i
\langle t_i^2,X_i \rangle \biggr| & \leq & c_1(L,\delta) \sqrt{\frac{m}{k}
\log({ek}/{m})}\E_g \sup_{t \in T} \Biggl|\sum_{i=1}^n g_i
t_i^2 \Biggr|
\\
& \leq & c_1(L,\delta) \sqrt{\frac{m}{k}
\log({ek}/{m})}\E_g \sup_{t \in T} \Biggl|\sum_{i=1}^n g_i
t_i \Biggr|.
\end{eqnarray*}
Also, recall that $P_I^\sigma F_1^{m} \subset c_2(L,\delta)
\gamma_{2,\tau_m}(F,L_2) W_{m}$, and, just as in the proof of Theorem
\ref{thm:herd} and by the isotropicity of $\mu$,
\[
\gamma_2(P_I^\sigma F_1^{m},| \cdot|) \leq c_3
\gamma_{2,\tau_m}(F,L_2) \sqrt{m} \leq c_3 \gamma_2(F,L_2)
\sqrt{m}=c_3 \gamma_2(T, | \cdot| )\sqrt{m}.
\]
Applying the Majorizing Measures theorem,
\[
\gamma_2(T,| \cdot|) \leq c_4\E\sup_{t \in T}
\Biggl|\sum_{i=1}^n g_i t_i \Biggr|,
\]
and thus
\begin{eqnarray*}
\frac{1}{\sqrt{k}} \E\sup_{t \in T} \biggl|\sum_{i \in I} g_i
\langle t_i^1,X_i \rangle \biggr| & \leq & \frac{c_4}{\sqrt{k}}
\gamma_2(P_I^\sigma F_1^{m},| \cdot|)
\\
& \leq & c_5 \sqrt{\frac{m}{k}} \E\sup_{t \in T} \Biggl|\sum_{i=1}^n
g_i t_i \Biggr|,
\end{eqnarray*}
as claimed.
\end{pf*}

To put Corollary \ref{cor:Rn} in the right context, even if one
considers the case where $\mu$ is the canonical Gaussian measure on
$\R^n$, the standard concentration estimate for the norm of a
Gaussian vector around its mean (used in \cite{Mi} to prove the
result for $I=\{1,\ldots,k\}$) is not strong enough to allow a uniform
control over all subsets of $\{1,\ldots,k\}$. What allows one to bypass
this obstacle and obtain a result even in a subgaussian setup (in
which case such a concentration result does not exist, and thus,
even the result for $I=\{1,\ldots,k\}$ is not obvious) is the
application of a cardinality-sensitive deviation argument rather
than a concentration based method.

Note that the logarithmic term in Corollary \ref{cor:Rn} cannot be
removed. For example, if $T=\{t\}$ and $\mu$ is the canonical
Gaussian measure on $\R^n$ then the vector $(\langle X_i,t \rangle)_{i=1}^k$
has the same distribution as $|t|(\bar{g}_i)_{i=1}^k$, where\vspace*{1pt}
$(\bar{g}_i)_{i=1}^k$ are independent standard Gaussian variables
[that are also independent of $(g_i)_{i=1}^n$]. Recall that $E_m$ is
the collection of subsets of $\{1,\ldots,k\}$ of cardinality $m$ and
observe that
\begin{eqnarray*}
\E_X \sup_{I \in E_m} \E_g \biggl|\sum_{i \in I} g_i \langle X_i,t
\rangle
\biggr| & \sim & \E_X \sup_{I \in E_m} \biggl(\sum_{i \in I}
\langle X_i,t \rangle^2 \biggr)^{1/2}
\\
& = &|t| \E\sup_{I \in E_m} \biggl(\sum_{i \in I} (\bar{g}_i)^2
\biggr)^{1/2} = |t| \E\Biggl(\sum_{i=1}^m
(\bar{g}_i^*)^2\Biggr)^{1/2}
\\
& \sim & |t| \sqrt{m \log({ek}/{m})},
\end{eqnarray*}
where the last assertion is the second part of Lemma
\ref{lemma:gauss-arr}. Therefore, with probability at least $c_1$,
there will be some $I \in E_m$ for which
\[
\frac{1}{\sqrt{k}}\E_g \biggl|\sum_{i \in I} g_i \langle X_i,t
\rangle
\biggr|
\geq c_2|t| \sqrt{ \frac{m}{k} \log({ek}/{m}) },
\]
showing that indeed, one cannot remove the logarithmic term.

\section{\texorpdfstring{Proof of Theorem \protect\ref{theoA}}{Proof of Theorem A.}} \label{sec:main}
As we explained in previous sections, our method of selecting signs
in a way that is better than choosing typical signs depends on
two properties. One is that the complexity of the set (as
captured, e.g., by $\gamma_{2,s}$ or the metric entropy of
the set) is small, and the other is that the set is in a good
position (e.g., if coordinate projections shrink the
set's complexity). Our results thus far indicate that for a
subgaussian class and a typical $\sigma=(X_i)_{i=1}^k$, $P_\sigma F$
is essentially a set in a good position. Thus, it seems likely that
the ability to choose signs that outperform the typical behavior of
signs will be governed solely by the complexity of $F$. As Theorem
\ref{theoA}, which we reformulate below, shows, this is indeed the case.

Although the proof of Theorem \ref{theoA} is rather technical, the basic idea
behind it is simple. It follows from a combination of the two main
results of the previous sections. First of all, that a typical
coordinate projection of a subgaussian class is contained in the
Minkowski sum of a small set and a set that satisfies a strong
shrinking property. Second, that the discrepancy of sets that
satisfy a shrinking property may be bounded in a nontrivial manner
using their metric complexity.
\begin{Theorem} \label{thm:main}
For any $0<\delta<1$, $0<\rho<1/2$ and $L>0$ there are constants
$c_1$ and $c_2$ that depend on $\delta$, $\rho$ and $L$ and for
which the following holds. Let $F \subset L_2(\mu)$ be an
$L$-subgaussian class, consisting of mean zero functions. Then, for
every $k$ there is a set $\mathcal{A}_k \subset\Omega^k$ with
$\mu^k(\mathcal{A}_k) \geq1-\delta$ such that for every $(X_1,\ldots,X_k)
\in\mathcal{A}_k$ and every $I \subset\{1,\ldots,k\}$,
\[
\inf_{(\eps_i)_{i \in I}} \sup_{f \in F} \biggl| \sum_{i \in I}
\eps_i f(X_i) \biggr| \leq\sqrt{|I|} a_{|I|},
\]
where for every $n \leq k$
\[
a_n \leq c_1\biggl(\gamma_{2,\log_2 \log_2 (c_2n)}(F,L_2) \cdot
\sqrt{\log(ek/n)}+\operatorname{diam}(F,L_2)
\frac{\log{k}}{n^{{1/2}-\rho}}\biggr).
\]
\end{Theorem}

Before proving Theorem \ref{thm:main}, let us recall the following
notation. For every integer $m$, $s_m$ is the largest integer $s$
such that $2^{2^{s+1}} \leq\ka_3 m$. If $m \leq k$, then $\tau_m$ is
the first integer for which $2^{2^s} \geq\exp(m \cdot\log(ek/m))$.
In particular, for every $1 \leq m \leq k$, $\tau_m \geq\log_2
\log_2 k$ (but of course, $\tau_m$ could be much larger). We will
also say that for $\sigma=(X_1,\ldots,X_k)$, a function class $F$
satisfies the shrinking property on $I \subset\{1,\ldots,k\}$ with a
constant $c$ if for every $f,h \in F$,
\[
\|f-h\|_{L_2^{I}} \leq c \sqrt{\log(ek/|I|)} \|f-h\|_{L_2}.
\]
\begin{pf*}{Proof of Theorem \ref{thm:main}}
Fix $0<\delta<1$ and consider $\sigma=(X_1,\ldots,X_k)$ for which the
assertions of Theorem \ref{thm:decomposition-of-F} hold. Fix any
integer $n \leq k$ and let $I_0 \subset\{1,\ldots,k\}$ be of cardinality
$n$. Using the notation of Theorem \ref{thm:decomposition-of-F}, we may
decompose $F \subset F_1^n + F_2^n$, where $P_{I_0}^\sigma F_1^n
\subset c_1 \gamma_{2,\tau_n}(F,L_2) W_n$, $F_2^n \subset F$, and
$F_2^n$ satisfies the shrinking property on every $I \subset
\{1,\ldots,k\}$ of cardinality $n$ with a constant $c=c(L,\delta)$---and
in particular, it does so on $I_0$.

For every $f \in F$, choose $f_1 \in F_1^n$ and $f_2 \in F_2^n$ such
that $f=f_1+f_2$. Hence, for every $(\eta_i)_{i \in I_0} \in
\{-1,0,1\}^{I_0}$ and any $f \in F$,
\[
\biggl| \sum_{i \in I_0} \eta_i f(X_i) \biggr| \leq
c_1\gamma_{2,\tau_n}(F,L_2) \sqrt{n} + \biggl| \sum_{i \in I_0}
\eta_i f_2(X_i) \biggr|.
\]

Let $(F_s)_{s=1}^\infty$ be an admissible sequence of $F_2^n$ which
will be specified later and set $\pi_s(f_2)$ to be a nearest point
to $f_2$ in $F_s$. As in Corollary
\ref{cor:choice-of-Q-for-gamma-2}, if
\[
Q_s^0 = \ka_4
\cases{
\exp(-\ka_5n^{1/2}), &\quad if $s<s_n$,\cr
1, &\quad if $s=s_n$,\cr
2^{s/2}, &\quad if $s>s_n$,}
\]
then there exist $(\eta_i^0)_{i \in I_0} \in\{-1,0,1\}^{I_0}$ such
that $n/4 \leq|\{i \dvtx \eta_i^0=0\}| \leq3n/4$ and for every
$(f_2(X_i))_{i \in I_0} \in P_{I_0}^\sigma F_2^n$,
\[
\biggl|\sum_{i \in I_0} \eta_i^0 f_2(X_i) \biggr| \leq\sum
_{s=1}^\infty
Q_s^0 \|\pi_s(f_2)-\pi_{s-1}(f_2)\|_{\ell_2^{I_0}}.
\]
Since functions in $F_2^n$ satisfy the shrinking property with a
constant $c$, then for every $f \in F_2^n$
\[
\|\pi_s(f)-\pi_{s-1}(f)\|_{\ell_2^{I_0}} \leq c\sqrt{n \log
(ek/n)} \|\pi_s(f)-\pi_{s-1}(f)\|_{L_2},
\]
implying that
\[
\biggl|\sum_{i \in I_0} \eta_i^0 f_2(X_i) \biggr| \leq c\sqrt{n
\log
(ek/n)} \sum_{s=1}^\infty Q_s^0
\|\pi_s(f_2)-\pi_{s-1}(f_2)\|_{L_2}.
\]
Let $I_1=\{i \in I_0\dvtx \eta_i=0\}$ and continue in the same manner:
first decompose $F \subset F_1^{|I_1|}+F_2^{|I_1|}$, then apply the
fact that
$P_{I_1}^\sigma F_1^{|I_1|}$ is contained in an appropriate weak
$\ell_2$ ball, and finally, since $F_2^{|I_1|}$ satisfies the
shrinking property on $I_1$, use
Corollary~\ref{cor:choice-of-Q-for-gamma-2} again, and so on.

As a result of iterating this argument, there are nested subsets of
$\{1,\ldots,k\}$, $(I_j)_{j=0}^{j_0}$, with $|I_0|=n$ and of
cardinalities
\[
\frac{|I_{j}|}{4} \leq|I_{j+1}| \leq\frac{3}{4}|I_{j}|,\qquad
1\leq|I_{j_0}| \leq10,
\]
and vectors $(\eta_i^j)_{i \in I_j} \in\{-1,0,1\}^{I_{j}}$, such
that $I_{j+1}=\{i \dvtx \eta_i^{j}=0\}$ with the following property. For
every $0 \leq j \leq j_0$, let
\[
Q_s^j = \ka_4
\cases{
\exp(-\ka_5|I_j|^{1/2}), &\quad if $s<s_{|I_j|}$,\cr
1, &\quad if $s=s_{|I_j|}$,\cr
2^{s/2}, &\quad if $s>s_{|I_j|}$,}
\]
and for
every $f \in F$, $f=f_1^j+f_2^j$, $f_1^j \in F_1^{|I_j|}$, $f_2^j \in
F_2^{|I_j|}$ one has
\begin{eqnarray*}
\biggl|\sum_{i \in I_j} \eta_i^{j} f(X_j) \biggr| &\leq&
\sum_{s=1}^\infty Q_s^j
\|\pi_s(f_2^j)-\pi_{s-1}(f_2^j)\|_{\ell_2^{I_j}} + c_1\sqrt{|I_j|}
\gamma_{2,\tau_{|I_j|}}(F,L_2)
\\
&\leq& c\sqrt{|I_j| \log(ek/|I_j|)} \sum_{s=1}^\infty
Q_s^j \|\pi_s(f_2^j)-\pi_{s-1}(f_2^j)\|_{L_2} \\
&&{} +
c_1\sqrt{|I_j|}\gamma_{2,\tau_{|I_j|}}(F,L_2).
\end{eqnarray*}
Therefore, there are signs $(\eps_i)_{i=1}^n \in\{-1,1\}^n$ such
that,
%
\begin{eqnarray} \label{eq:disc-error-subgauss-proj}
&&
\sup_{f \in F} \Biggl|\sum_{i=1}^n \eps_i f(X_i) \Biggr|
\nonumber
\\
&&\qquad\leq c\sum_{j=0}^{j_0} \sqrt{|I_j| \log
(ek/|I_j|)} \sup_{f \in F} \sum_{s=1}^\infty Q_s^j
\|\pi_s(f)-\pi_{s-1}(f)\|_{L_2}
\\
&&\qquad\quad{} + c_1\sum_{j=0}^{j_0} \sqrt{|I_j|}\gamma_{2,\tau_{|I_j|}}(F,L_2) +
c_2\operatorname{diam}(F,L_2)\log(ek),\nonumber
\end{eqnarray}
where the last term comes from a trivial estimate on the discrepancy
of a projection of $F$ onto the set of coordinates $\{i \in I_{j_0}
\dvtx \eta_i^{j_0} = 0\}$ and the shrinking phenomenon.

To complete the proof, one has to bound
(\ref{eq:disc-error-subgauss-proj}) from above. To that end, set
$b_j=|I_j|$ and recall that $(1/4)^j n \leq b_j \leq(3/4)^j n$. To
estimate the second term in (\ref{eq:disc-error-subgauss-proj}),
since $b_j \geq1$ then $\tau_{b_j} \geq\log_2 \log_2 k$.
Therefore,
\[
\sum_{j=0}^{j_0} \sqrt{|I_j|}\gamma_{2,\tau_{|I_j|}}(F,L_2) \leq c_3
\sqrt{n} \gamma_{2,\log_2 \log_2 k}(F,L_2).
\]

Turning our attention to the first term in
(\ref{eq:disc-error-subgauss-proj}), for every $1 \leq\ell\leq
s_n$ let $U_\ell=\{j\dvtx s_{b_j} =\ell\}$, set $b^+_\ell= \max\{ b_j \dvtx j \in U_\ell\}$ and $b^-_\ell=\min\{ b_j \dvtx j \in U_\ell\}$. In other
words, $U_\ell$ consists of all the integers $j$ for
which $s_{|I_j|}=s_{b_j}=\ell$; $b^+_\ell$ is the largest cardinality
of such
a set and $b_\ell^-$ is the smallest one. Since
\[
\ka_3^{-1} 2^{2^{\ell+1}} \leq b_\ell^{-} \leq b_\ell^+ \leq
\min\{\ka_3^{-1} 2^{2^{\ell+2}},n\},
\]
then for every $j \in U_\ell$, the sequence $(Q_s^j)_{s>0}$
satisfies
\[
Q_s^j \leq\ka_4
\cases{
\exp(-\ka_5 \ka_3^{-1/2} \cdot2^{2^{\ell}}),
&\quad if $s<\ell$,\cr
1, &\quad if $s=\ell$, \cr
2^{s/2}, &\quad if $s>\ell$,}
\]
and we denote this sequence $(Q_s^\ell)_{s>0}$. Since $b_j$ decays
exponentially, then
\begin{eqnarray*}
&& \sum_{j=0}^{j_0} \sqrt{|I_j| \log(ek/|I_j|)} \sup_{f
\in F} \sum_{s=1}^\infty Q_s^j \|\pi_s(f)-\pi_{s-1}(f)\|_{L_2}
\\
&&\qquad = \sum_{\ell=1}^{s_n} \sum_{b_j \in U_\ell} \sqrt{b_j \log
(ek/b_j)} \sup_{f \in F} \sum_{s=1}^\infty Q_s^j
\|\pi_s(f)-\pi_{s-1}(f)\|_{L_2}
\\
&&\qquad \leq c_4 \sum_{\ell=1}^{s_n} \sqrt{b^+_\ell\log
(ek/b^+_\ell
)} \sup_{f \in F} \sum_{s=1}^\infty Q_s^\ell
\|\pi_s(f)-\pi_{s-1}(f)\|_{L_2}.
\end{eqnarray*}

Set $d_\ell= \sqrt{b^+_\ell\log(ek/b^+_\ell)}$, fix
$0< \rho< 1/2$ and let $\ell_1$ be the largest integer such that
$\ka_3^{-1}2^{2^{\ell_1+2}} \leq n^{2\rho}$. Then, for every $\ell
\leq\ell_1$, $b_\ell^+ \leq\ka_3^{-1}2^{2^{\ell_1+2}} \leq
n^{2\rho}$ and for $\ell> \ell_1$, $b_\ell^+ \leq n$. Observe that
for every $s,\ell$, $Q_s^\ell\leq\ka_4 2^{s/2}$ and for every $f
\in F$, $\|\pi_s(f)-\pi_{s-1}(f)\|_{L_2} \leq2\operatorname{diam}(F,L_2)$.
Therefore,
\begin{eqnarray*}
&& \sum_{\ell=1}^{s_n} d_\ell\cdot\sup_{f \in F} \sum_{s=1}^\infty
Q_s^\ell\|\pi_s(f)-\pi_{s-1}(f)\|_{L_2}
\\
&&\qquad\leq \sum_{\ell=1}^{s_n} d_\ell\cdot\sup_{f \in F}
\sum_{s=1}^{\ell_1} Q_s^\ell\|\pi_s(f)-\pi_{s-1}(f)\|_{L_2}\\
&&\qquad\quad{} +
\sum_{\ell=1}^{s_n} d_\ell\cdot\sup_{f \in F}
\sum_{s=\ell_1+1}^\infty Q_s^\ell\|\pi_s(f)-\pi_{s-1}(f)\|_{L_2}
\\
&&\qquad\leq 2\operatorname{diam}(F,L_2) \sum_{\ell=1}^{s_n} \sum_{s=1}^{\ell_1}
d_\ell Q_s^\ell\\
&&\qquad\quad{} + \ka_4 \sum_{\ell=1}^{s_n} d_\ell\cdot\sup_{f\in
F} \sum_{s=\ell_1+1}^\infty2^{s/2} \|\pi_s(f)-\pi_{s-1}(f)\|_{L_2}
\\
&&\qquad\leq 2\operatorname{diam}(F,L_2) \sum_{\ell=1}^{s_n} \sum_{s=1}^{\ell_1}
d_\ell Q_s^\ell+ c_5(\rho) \sqrt{n \log{(ek/n)}} \cdot
\gamma_{2,\ell_1}(F,L_2)
\end{eqnarray*}
for an almost optimal choice of $(F_s)_{s=\ell_1}^\infty$.

Now, for every $s \leq\ell_1$ and using that $b_\ell^+ \leq
n^{2\rho}$ for $\ell< \ell_1$ and $b_\ell^+ \leq n$ for $\ell\geq
\ell_1$, it is evident that
\begin{eqnarray*}
\sum_{\ell=1}^{s_n} d_\ell Q_s^\ell & = & \sum_{\ell=1}^{s_n}
\sqrt{b^+_\ell\log(ek/b^+_\ell)} Q_s^\ell
\\
& \leq & n^\rho\sqrt{\log({ek}/{n^{2\rho}})} \sum_{\ell<
\ell_1} Q_s^\ell+ \sqrt{n\log({ek}/{n})} \sum_{\ell
\geq
\ell_1} Q_s^\ell= (*).
\end{eqnarray*}
Note that $2^{\ell_1+3} \geq2\rho\log_2{c_6n}$ and thus
$2^{\ell_1+1} \geq(\rho/2)\log_2{c_6n}$. Therefore, there is an
absolute constant $c_7$ such that
if $s \leq\ell_1$ then
\begin{eqnarray*}
\sum_{\ell\leq\ell_1} Q_s^\ell & = & \sum_{\ell\leq s} Q_s^\ell+
\sum_{\ell= s+1}^{\ell_1} Q_s^\ell
\\
& \leq &c_7 \bigl(s2^{s/2} + \exp(-c_3 2^{2^{s}}) \bigr) \leq
2c_7 s2^{s/2}
\end{eqnarray*}
and
\[
\sum_{\ell> \ell_1} Q_s^\ell\leq c_7 \exp(-c_3 2^{
2^{\ell_1}}) \leq c_7\exp(-c_8n^{\rho/2}).
\]
Hence, there is a constant $c_9(\rho)$ such that for every $s
\leq\ell_1$,
\[
(*) \leq c_9(\rho) s2^{s/2}n^\rho\sqrt{\log
({ek}/{n^{2\rho}})},
\]
and thus,
\begin{eqnarray*}
\operatorname{disc}(P_I^\sigma F) & \leq & c_2\operatorname{diam}(F,L_2)\log(ek)+ c_3
\sqrt{n} \gamma_{2,\log_2 \log_2 k}(F,L_2)
\\
&&{} + c_{10} \operatorname{diam}(F,L_2) \cdot\ell_12^{\ell_1/2} n^\rho\sqrt
{\log
({ek}/{n^{2\rho}})}
\\
&&{} + c_{10} \gamma_{2,\ell_1}(F,L_2) \cdot\sqrt{n \log
({ek}/{n})}.
\end{eqnarray*}
Since $(\rho/4) \log_2 c_6n \leq2^{\ell_1} \leq(\rho/2) \log_2
c_6n$, the claim follows.
\end{pf*}
\begin{Corollary}
Let $0<\rho<1/2$. Under the assumptions of Theorem \ref{thm:main}
and using its notation, for every $\sigma\in\mathcal{A}_k$
\[
\operatorname{disc}(P_\sigma F) \leq c_1 \bigl(\sqrt{k} \cdot
\gamma_{2,\log_2 \log_2 (c_2k)}(F,L_2) + k^\rho\operatorname{diam}(F,L_2)\bigr)
\]
and
\[
\operatorname{Hdisc}(P_\sigma F)= \sup_{I \subset\{1,\ldots,k\}}
\inf_{(\eps_i)_{i=1}^k} \sup_{f \in F} \biggl| \sum_{i \in I} \eps_i
f(X_i) \biggr| \leq\sup_{1 \leq n \leq k} a_n \sqrt{n}.
\]
In particular, if $\lim_{s \to\infty}
\gamma_{2,s}(F,L_2) =0$ (i.e., if $F$ is $\mu$-pregaussian), then
\[
\frac{1}{\sqrt{k}} \operatorname{Hdisc}\bigl(\{
(f(X_i))_{i=1}^k \dvtx f \in F \}\bigr)
\]
converges in probability to 0.
\end{Corollary}

Let us mention\vspace*{1pt} once again that the reason that Theorem
\ref{thm:main} is meaningful is because for a typical
$(X_i)_{i=1}^k$, a class of mean zero functions that is
$L$-subgaussian satisfies that
\[
c_1 \sigma_F \sqrt{k} \leq\E_\eps\sup_{f \in F} \Biggl|\sum_{i=1}^k
\eps_i f(X_i)\Biggr| \leq c_2(L) \gamma_2(F,L_2) \sqrt{k}.
\]
Thus, there is a true gap between the discrepancy (or even the
hereditary discrepancy) of a typical coordinate projection and the
average over signs of a coordinate projection of a pregaussian,
subgaussian class $F$.

\section{Equivalence for large sets} \label{sec:upp-low}
In this section, our aim is to show that if $F$ is a subgaussian
class that indexes a bounded Gaussian process, then the reason for
the gap between the expectation over signs of a random coordinate
projection and the infimum over signs is indeed that $\lim_{s \to
\infty}
\gamma_{2,s}(F,L_2) =0$.

To be more precise, we show the following.
\begin{Theorem} \label{thm:optimal}
For every $0<\delta<1$ and $A,B,L>0$ there is a constant
$c(\delta,A,B,L)$ for which the following holds. Let $F \subset
B(L_2(\mu))$ be a class of mean
zero functions such that $\operatorname{absconv}(F)$ is $L$-subgaussian.
If $\gamma_2(F,L_2(\mu)) \leq A < \infty$ and
if the entropy numbers satisfy that
\[
\limsup_{j \to\infty} j^{1/2}e_j(\operatorname{absconv}(F),L_2(\mu))=
B>0,
\]
then there is a sequence of integers $(k_i)_{i=1}^\infty$ tending to
infinity, such that for every $i$, with probability at least
$1-\delta$ in $\Omega^{k_i}$,
\[
\operatorname{Hdisc}(P_\sigma F) \geq c(\delta,A,B,L)\sqrt{k_i},
\]
where $\sigma=(X_1,X_2,\ldots,X_{k_i}) \in\Omega^{k_i}$ is selected
according to $\mu^{k_i}$. In particular, $\operatorname{Hdisc}(P_\sigma F)/\sqrt{k}$ does not converge to $0$ in
probability.
\end{Theorem}

Observe that this is almost the reverse direction of Theorem \ref{theoA}.
Indeed, it is well known (see, e.g., \cite{Dud-book}, Chapter 9) that
there is
no entropic characterization of classes that index a bounded
Gaussian process which is not continuous; such a characterization is
given by a majorizing measures argument \cite{Tal87}. However,
because $\{G_f \dvtx f \in F\}$ is a bounded process with a covariance
structure endowed by $L_2(\mu)$, then by Sudakov's inequality (see,
e.g., \cite{LT}, Chapter 3),
\[
\log N(\eps,\operatorname{absconv}(F),L_2(\mu)) \leq c_1\biggl(\frac{\E\sup
_{f \in F}
G_f}{\eps}\biggr)^2.
\]
On the other hand, since $F \subset B(L_2(\mu))$ is not
$\mu$-pregaussian, one can show that
\[
\int_0^1 \sqrt{\log N(\eps,F,L_2(\mu))} \,d\eps= \infty.
\]

Thus, up to a logarithmic factor, the entropy numbers of $F$ are as
in Theorem \ref{thm:optimal}. Whether Theorem \ref{thm:optimal}
remains true using only the assumption that $\limsup_{s \to\infty}
\gamma_{2,s}(F,L_2(\mu)) >0$ is not clear.

The idea behind the proof of Theorem \ref{thm:optimal} is to find a
cube in a typical coordinate projection of $\operatorname{absconv}(F)$. We will
first show
that if $\operatorname{absconv}(F)$ has a ``large'' separated set with respect to the
$L_2(\mu)$ metric at scale \mbox{$\sim$}$1/\sqrt{k}$, then its typical
coordinate projection of dimension $k$ contains a cubic structure of
dimension \mbox{$ \sim$}$k$ and scale \mbox{$\sim$}$1/\sqrt{k}$. The cubic structure
we will
be interested in is captured by the combinatorial dimension.
\begin{Definition} \label{def:comb-dim}
Let $F$ be a class of functions on $\Omega$. For every $\eps>0$,
a~set $\sigma=\{ x_1,\ldots,x_j \} \subset\Omega$ is said to be
$\eps$-shattered by $F$ if there is some function $s\dvtx\sigma\to\R$,
such that for every $I \subset\{1,\ldots,j\}$ there is some $f_I \in
F$ for which $f_I(x_i) \geq s(x_i)+\eps$ if $i \in I$, and $f_I(x_i)
\leq s(x_i)-\eps$ if $i \notin I$. Define the combinatorial
dimension at scale $\eps$ by
\[
\operatorname{VC}(F,\eps)=\sup\{|\sigma| |  \sigma
\subset
\Omega, \sigma \mbox{ is }  \eps\mbox{-shattered  by}  F
\}.
\]
\end{Definition}

Note that if $F$ is a $\{0,1\}$-class of functions then
$\operatorname{VC}(F)=\operatorname{VC}(F,1/2)$. Also, in a similar way
one may define the combinatorial dimension of a subset of $\R^n$, when
each vector is viewed as a function defined on $\{1,\ldots,n\}$.

It is standard to verify that if $\operatorname{VC}(F,\eps) \geq m$,
then the coordinate projection $P_\tau F$, defined by the shattered set
$\tau$, contains a subset of cardinality $\exp(cm)$ which is
$c_1\eps$-separated with respect to the $L_2(\mu_\tau)$ norm (recall
that $\mu_\tau$ is the uniform probability measure supported on
$\tau$), and that $\operatorname{disc}(P_\tau F) \geq c_2m\eps$ (see
Lem\-ma~\ref{lemma:disc-vs-VC}). As we mentioned in the
\hyperref[intro]{Introduction}, the reverse direction is also true, and
if $F \subset B(L_\infty(\Omega))$ contains a large well-separated set
in $L_2(\mu)$ that it must have a large combinatorial dimension at a
scale that is proportional to the scale of the separation (see
\cite{MenVer} for an exact statement and proof). A fact that will be
used here and which is based on this reverse direction is the
following.
\begin{Theorem}[\cite{MenVer}] \label{thm:MV}
There exist\vspace*{1pt} absolute constants $c_1$ and $c_2$ for which the
following holds. Let $V \subset B_\infty^k$ and assume that $\E
{\sup_{v \in V}} |{\sum_{i=1}^k \eps_i v_i} | \geq\delta k$. Then,
\[
\VC(V,c_1\delta) \geq c_2\delta^2k.
\]
\end{Theorem}

Hence, the only reason that ${\E
\sup_{v \in V} }|{\sum_{i=1}^k \eps_i v_i} |$ is almost extremal is
that $V$ contains a large cube in
a high-dimensional coordinate projection.

The key observation of this section is the following theorem.
\begin{Theorem} \label{thm:cube}
For every $A,B,L>0$ and $0<\delta<1$ there exist constants $c_1$ and
$c_2$ that depend on $A,B,L$ and $\delta$ for which the following
holds. Let $F \subset B(L_2(\mu))$ be a convex, symmetric,
$L$-subgaussian set of mean zero functions. Suppose that\vspace*{1pt}
$\gamma_2(F,\psi_2) \leq A<\infty$ and that there is some $k$ for
which $e_k(F,L_2(\mu)) \geq B/\sqrt{k}$. Then, there is a set
$\Sigma\subset\Omega^k$ such that $\mu^k(\Sigma) \geq1-\delta$ and
for every $\sigma\in\Sigma$,
\[
\VC\biggl(P_\sigma F, \frac{c_1}{\sqrt{k}}\biggr) \geq c_2k.
\]
\end{Theorem}

Theorem \ref{thm:cube} implies Theorem \ref{thm:optimal} because of
the next lemma.
\begin{Lemma} \label{lemma:disc-vs-VC}
If $T \subset\R^n$, then
\[
\operatorname{Hdisc}(T) \geq\sup_{\delta>0} \delta\operatorname{VC}(\operatorname{absconv}(T),\delta).
\]
\end{Lemma}
\begin{pf}
First, note that $\operatorname{Hdisc}(T)=\operatorname{Hdisc}(\operatorname{absconv}(T))$, and thus
we may assume that $T$ is convex and symmetric. Now,
let $I \subset\{1,\ldots,n\}$ be $\delta$-shattered by
$T$ with the level function $s$. Fix
$(\eps_i)_{i \in I} \in\{-1,1\}^{|I|}$ and without loss of
generality assume that $\sum_{i \in I} \eps_i s_i \geq0$. Since $I$ is
$\delta$-shattered by $T$, there is some $t^\prime\in T$ for which
$t_i^\prime\geq s_i + \delta$ when $\eps_i =1$ and $t_i^\prime\leq
s_i -\delta$ when $\eps_i=-1$. Thus,
\[
\sup_{t \in T} \biggl|\sum_{i \in I} \eps_i t_i\biggr| \geq
\biggl|\sum_{i \in I} \eps_i (t^\prime_i - s_i) + \sum_{i \in I}
\eps_i s_i \biggr| \geq\biggl|\sum_{i \in I} \eps_i (t^\prime_i -
s_i) \biggr| \geq|I| \delta,
\]
as claimed.
\end{pf}

Hence, from here on we may assume without loss of generality that the
class $F$ is convex and symmetric, and that it is $L$-subgaussian.

The proof of Theorem \ref{thm:cube} requires several additional
facts. To formulate them, denote for $V \subset\R^n$
\[
\ell_*(V) = \E\sup_{v \in V} \Biggl|\sum_{i=1}^n g_i v_i \Biggr|,
\]
and if $A,B \subset\R^n$, set $N(A,B)$ to be the minimal number of
translates of $B$ needed to cover $A$.

The first lemma we need is taken from \cite{LMPT}.
\begin{Lemma} \label{lemma:LMPT}
Let $V \subset\R^k$ be a convex, symmetric set. For $\rho>0$, set
$V_\rho= V \cap\rho B_2^k$ and $F(\rho)=\ell_*(V)/\ell_*(V_\rho)$.
Then,
\[
N(V,8\rho B_2^k) \leq
\exp\biggl(2\biggl(\frac{\ell_*(V_\rho)}{\rho}\biggr)^2\log
(6F(\rho))\biggr).
\]
\end{Lemma}

The second result was proved in \cite{MPT} (Theorem 2.3). Although
it was formulated there for subsets of $\R^n$, its proof shows that
the claim is true for any subgaussian class of functions. It implies
that a
random coordinate projection of $F$, viewed as a mapping between
$L_2(\mu)$ and $L_2^k$, is almost norm preserving for functions with
a sufficiently large $L_2(\mu)$ norm.
\begin{Theorem} \label{thm:MPT}
There exist absolute constants $c_1$ and $c_2$ for which the
following holds. Let $F \subset L_2(\mu)$ be a convex, symmetric,
$L$-subgaussian class
of functions. For every $\theta>0$ and any positive integer $k$, set
\[
r_k(\theta)=\inf\biggl\{\rho\dvtx \rho\geq\frac{\gamma_2(F \cap\rho
B(L_2(\mu)),\psi_2)}{\theta\sqrt{k}} \biggr\}.
\]
Then, with probability at least $1-2\exp(-c_1\theta^2k/L^4)$, for
every $f \in F$ such that $\|f\|_{L_2(\mu)} \geq
r_k(\theta/c_2L^2)$,
\[
(1-\theta)^{1/2} \|f\|_{L_2(\mu)} \leq\|f\|_{L_2^k} \leq
(1+\theta)^{1/2} \|f\|_{L_2(\mu)}.
\]
\end{Theorem}
\begin{Corollary} \label{cor:iso}
For every $L>0$, there are constants $\ka_{6}$ and $\ka_{7}$ that
depend only on $L$, for which the following holds. Let $F$ be an
$L$-subgaussian, convex and symmetric class of functions for which
$\gamma_2(F,\psi_2) \leq A < \infty$. Then, with probability at
least $1-2\exp(-\ka_{6}k)$, if $f \in F$ and $\|f\|_{L_2(\mu)} \geq
\ka_{7}A/\sqrt{k}$ then
\[
\sqrt{\tfrac1 2} \|f\|_{L_2(\mu)} \leq\|f\|_{L_2^k}
\leq\sqrt{\tfrac3 2} \|f\|_{L_2(\mu)}.
\]
In particular, if $H \subset F$ is an $\eps$-separated set in
$L_2(\mu)$ for $\eps>2\ka_{7}A/\sqrt{k}$ then with probability at
least $1-2\exp(-\ka_{6}k)$, $P_\sigma H$ is $\eps/4$-separated in
$L_2^k$.
\end{Corollary}
\begin{pf}
Let $c_1$ and $c_2$ be as in Theorem \ref{thm:MPT}. Observe that
\[
\gamma_2\bigl(F\cap\rho B(L_2(\mu)),\psi_2\bigr) \leq
\gamma_2(F,\psi_2) \leq A
\]
and apply Theorem \ref{thm:MPT} for $\theta=1/2$. Thus,
$r_k(\theta/c_2L) \leq c_3(L)A/\sqrt{k}$, implying that if
$c_4(L)=c_1/4L^4$ then with probability at least
$1-2\exp(-c_4(L)k)$, if $\|f\|_{L_2(\mu)} \geq c_3(L)A/\sqrt{k}$
then
\[
\tfrac{1}{2} \|f\|_{L_2(\mu)}^2 \leq\|f\|_{L_2^k}^2
\leq\tfrac{3}{2} \|f\|_{L_2(\mu)}^2.
\]
Turning to the second part, note that if $H \subset F$ is
$\eps$-separated in $L_2(\mu)$ for $\eps>2c_1(L)A/\sqrt{k}$, then
for every $h_1,h_2 \in H$, $f=(h_1-h_2)/2 \in F$ and
$\|f\|_{L_2(\mu)} \geq c_1(L)A/\sqrt{k}$. Thus, the second part
follows from the first one.
\end{pf}

Now we can formulate the first localization result, showing that the
richness of a typical coordinate projection comes from the
intersection of $F$ with a ball of radius $\sim1/\sqrt{k}$.
\begin{Theorem} \label{thm:intersection}
For every positive $A$, $B$, $L$ and $0<\delta<1$, there are
constants $c>1$, $c_1$ $c_2$ and $c_3$ depending on $A$, $B$, $L$
and $\delta$ for which the following holds. Let $F \subset
B(L_2(\mu))$ be a convex, symmetric, $L$-subgaussian class of mean
zero functions such that $\gamma_2(F,\psi_2) \leq A <\infty$. Fix an
integer $k$ and assume that $e_k(F,L_2(\mu)) \geq B/\sqrt{k}$. Then,
with probability at least $1-\delta-2\exp(-c_{1}k)$,
\[
\E_g \sup_{f \in F \cap{c_{2}}/{\sqrt{k}}B(L_2(\mu))}
\Biggl|\sum_{i=1}^{ck} g_i f(X_i) \Biggr| \geq c_{3}\sqrt{k}.
\]
\end{Theorem}
\begin{pf}
Since $F$ is $L$-subgaussian and by applying Sudakov's inequality,
we may assume without loss of generality that $A/B>1$.
Let $H$ be a maximal $B/\sqrt{k}$ separated set in $F$ with
$\log|H| \geq k$. Let $k^\prime=c^2k$ for a constant $c > 1$ to be
named later. Since $H$ is
$\eps=cB/\sqrt{k^\prime}$ separated in $L_2(\mu)$, then by Corollary
\ref{cor:iso}, with probability at least
$1-2\exp(-\ka_{6}k^\prime)=1-2\exp(-c_1(L)k)$, if $\eps\geq
2\ka_{7}A/\sqrt{k^\prime}$ then $P_\sigma H$ is $\eps/4$-separated in
$L_2^{k^\prime}$. Moreover, if $f \in F$ satisfies $\|f\|_{L_2(\mu)}
\geq
\ka_{7}A/\sqrt{k^\prime}$ then
%
\begin{equation} \label{eq:iso}
\tfrac{1}{2} \|f\|_{L_2(\mu)}^2 \leq\|f\|_{L_2^{k^\prime}}^2 \leq
\tfrac{3}{2} \|f\|_{L_2(\mu)}^2.
\end{equation}
Clearly, the condition on $\eps$ holds if $c \sim_L A/B$, and since $c>1$
it follows that $k^\prime> k$.

Consider the set $U =\operatorname{absconv}(H)$. By the Majorizing Measures
theorem and a simple application of Theorem
\ref{thm:decomposition-of-F}, with probability at least $1-\delta$
for $|\sigma|=k^\prime$,
%
\begin{equation} \label{eq:4.2}
\ell_*(P_\sigma U) \leq c_2 \gamma_2(P_\sigma U, | \cdot|) \leq
c_3(L,\delta) A \sqrt{k^\prime}.
\end{equation}
Let $\sigma=(X_i)_{i=1}^{k^\prime}$ be in the intersection of the
two events given by (\ref{eq:iso}) and (\ref{eq:4.2}), set $V =
P_\sigma U$ and note that
$B(L_2^{\sigma})=\sqrt{k^\prime}B_2^{k^\prime}$. Therefore,
%
\begin{eqnarray} \label{eq:4.3}
k &\leq&\log N\biggl(V,
\frac{\eps}{4}\sqrt{k^\prime} B_2^{k^\prime}\biggr) =\log N(V, c_4
B_2^{k^\prime})
\nonumber\\[-8pt]\\[-8pt]
&\equiv&\log N(V, 8\rho B_2^{k^\prime}),\nonumber
\end{eqnarray}
where $c_4 \sim_L A$ (and thus $\rho\sim_L A$ as well). If
$V_\rho=V \cap\rho B_2^{k^\prime}$, then by Lemma \ref{lemma:LMPT},
(\ref{eq:4.2}) and (\ref{eq:4.3}),
\[
k \leq2\biggl(\frac{\ell_*(V_\rho)}{\rho}\biggr)^2 \log(6F(\rho))
\leq2\biggl(\frac{\ell_*(V_\rho)}{\rho}\biggr)^2 \log
\biggl(\frac{c_3 A \sqrt{k^\prime}}{\ell_*(V_\rho)}\biggr).
\]

Solving this inequality for $\ell_*(V_\rho)$, it is evident that
there exists a constant $c_5 \sim_{L,\delta}
B/\sqrt{\log(c_6A^2/B^2)}$ (where $c_6$ depends on $L$ and
$\delta$) for which $\ell_*(V_\rho) \geq c_5 \sqrt{k^\prime}$. Since
$F$ is convex and symmetric and $H \subset F$, then
\[
V_\rho= P_\sigma(\operatorname{absconv}(H)) \cap\rho
B_2^{k^\prime} \subset P_\sigma\biggl(\biggl\{ f \in F \dvtx \|f\|_{L_2^{k^\prime}} \leq\frac{\rho}{\sqrt{k^\prime}}
\biggr\}\biggr)
\]
and by (\ref{eq:iso}),
\[
\biggl\{f \in F \dvtx \|f\|_{L_2^{k^\prime}} \leq
\frac{\rho}{\sqrt{k^\prime}} \biggr\} \subset\biggl\{ f \in F \dvtx \|f\|_{L_2(\mu)} \leq2 \frac{\max\{\rho,
\ka_{7}A\}}{\sqrt{k^\prime}} \biggr\}.
\]
Hence, there is a constant $c_7 \sim_{L,\delta} A$ for which with
probability at least $1-\delta-2\exp(-c_1k)$,
\[
V_\rho\subset P_\sigma\biggl(\biggl\{ f \in F \dvtx \|f\|_{L_2(\mu)}
\leq
\frac{c_7}{\sqrt{k^\prime}} \biggr\}\biggr),
\]
implying that
\[
c_5\sqrt{k^\prime} \leq\E_g \sup_{\{f \in F \dvtx \|f\|_{L_2(\mu)}
\leq
{c_7}/{\sqrt{k^\prime}}\}} \Biggl| \sum_{i=1}^{k^\prime} g_i
{f(X_i)} \Biggr|.
\]
\upqed\end{pf}

The next step in the proof of Theorem \ref{thm:cube} is a second
localization argument. Theorem \ref{thm:intersection} shows that under our
assumptions, there is a small ball (of radius \mbox{$\sim$}$1/\sqrt{k}$) in
$F$ that causes coordinate projections of $F$ of dimension $k$ to be
``rich.'' Now, one has to localize even further by truncating the
functions in $F_1=F \cap(c/\sqrt{k}) B(L_2(\mu))$.
\begin{Definition} \label{def:beta}
For every $\beta>0$ and every $f \in F$, let
\[
f_\beta^-= f \IND_{\{|f| \leq\beta\}}+\operatorname{sgn}(f) \beta
\IND_{\{|f| \geq\beta\}}
\]
and $f_\beta^+=f-f_\beta$. For every $\sigma=(X_1,\ldots,X_k)$ let
\[
V^-_\beta=\{(f_\beta^-(X_i))_{i=1}^k  \dvtx  f \in F\},\qquad
V^+_\beta=\{(f_\beta^+(X_i))_{i=1}^k  \dvtx  f \in F\}.
\]
\end{Definition}
\begin{pf*}{Proof of Theorem \ref{thm:cube}} First, by Theorem
\ref{thm:intersection}, with probability at least
$1-\delta-2\exp(-c_{1}k)$,
\[
\E_g \sup_{f \in F \cap{c_{2}}/{\sqrt{k}}B(L_2(\mu))}
\Biggl|\sum_{i=1}^{ck} g_i f(X_i) \Biggr| \geq c_{3}\sqrt{k},
\]
where $c>1$. Set
\[
H= F \cap\frac{c_2}{\sqrt{k}}B(L_2(\mu))
\]
and note that by the proof of Theorem \ref{thm:decomposition-of-F}
for the class $H$ and $m=ck$, each $h \in H$ can be written as
$h=h_1+h_2$, where $h_1 \in H - H \subset2H$ (by the convexity and
symmetry of $F$), and $h_2 \in H$. Moreover, if we write $H \subset
H_1 + H_2$ then with $\mu^{ck}$ probability $1-\delta$, $P_\sigma
H_1 \subset c_4 \gamma_{2,\tau_{ck}}(F,L_2)W_{ck} \subset c_4 A
W_{ck}$, where $c_4=c_4(L,\delta)$. By a standard concentration
argument---similar to the one used in Theorem \ref{thm:decomposition-of-F},
since $|H_2| \leq\exp(c_5k)$ then for every $h_2 \in H_2$,
$\|h_2\|_{L_2^{ck}} \leq c_6\|h_2\|_{L_2} \leq c_7/\sqrt{k}$. Thus,
$P_\sigma H_2 \subset(c_7/\sqrt{k})B_2^{ck}$, and since $B_2^{ck}
\subset W_{ck}$, then
\[
P_\sigma H \subset P_\sigma H_1 + P_\sigma H_2 \subset c_8 W_{ck},
\]
where $c_8=c_8(A,L,\delta)$.

Let $\sigma=(X_i)_{i=1}^{ck}$ for which the above estimates hold,
fix $\beta$ to be named later and let $V^+_\beta$ and $V^-_\beta$ be
as in Definition \ref{def:beta} for the set $H$. Consider the set
\[
W_k^\beta= \bigl\{x \in\R^{ck} \dvtx x_i^* \leq\bigl(c_8/\sqrt{i}\bigr) - \beta
\mbox{ for }  i \leq(c_8/\beta)^2,  x_i^* = 0  \mbox{ for }  i >
(c_8/\beta)^2 \bigr\}
\]
and observe that $V^+_\beta\subset W_{ck}^\beta$. Therefore, if we
set $i_\beta=(c_8/\beta)^2$ and select $\beta$ to satisfy that $1
\leq i_\beta\leq ck$ then,
\begin{eqnarray*}
\E\sup_{v \in V^+_\beta} \Biggl|\sum_{i=1}^{ck} g_i v_i \Biggr|
&\leq& \E\sup_{w \in W_{ck}^\beta} \Biggl|\sum_{i=1}^{ck} g_i v_i
\Biggr| \leq\E\sum_{i=1}^{i_\beta} \frac{c_8}{\sqrt{i}} g_i^*
\\
&\leq& c_9 \sqrt{i_\beta\log(ek/i_\beta)} \leq
\frac{c_{3}}{2}\sqrt{k}
\end{eqnarray*}
for an appropriate choice of $\beta\sim c_3/\sqrt{k}$. Since
$V=P_\sigma H \subset V_\beta^+ + V_\beta^-$, then
\begin{eqnarray*}
c_3\sqrt{k} & \leq & \E\sup_{v \in V} \Biggl|\sum_{i=1}^{ck} g_i v_i
\Biggr| \\
&\leq& \E\sup_{v \in V^-_\beta} \Biggl|\sum_{i=1}^{ck} g_i v_i
\Biggr| +\E\sup_{v \in V^+_\beta} \Biggl|\sum_{i=1}^{ck} g_i v_i\Biggr|
\\
& \leq& \E\sup_{v \in V^-_\beta} \Biggl|\sum_{i=1}^{ck} g_i v_i
\Biggr| + \frac{c_3}{2} \sqrt{k}.
\end{eqnarray*}
Therefore,
\[
\E\sup_{v \in\beta^{-1}V^-_\beta} \Biggl|\sum_{i=1}^{ck} g_i v_i
\Biggr| \geq\frac{c_3\sqrt{k}}{2 \beta} \geq c_{10} k.
\]
Note that
\[
\beta^{-1}V^-_\beta= \biggl\{ \sum_{\{i\dvtx |f(X_i)| \leq\beta\}}
\beta^{-1}f(X_i)e_i + \sum_{\{i\dvtx |f(X_i)| > \beta\}}
\sgn(f(X_i))e_i \dvtx f \in H\biggr\} \subset B_\infty^{ck}.
\]
Therefore, by the optimal estimate in the sign-embedding theorem \cite{MenVer},
there are constants $c_{11} \sim c_{10}^2$ and $c_{12} \sim c_{10}$
such that
\[
\VC(\beta^{-1}V^-_\beta, c_{11}) \geq c_{12}k.
\]
In other words, there is a set $I \subset\{1,\ldots,ck\}$, $|I| \geq
c_{12}k$ and a vector $(s_i)_{i \in I}$ such that for every $J
\subset I$, there is $v_J \in V^-_\beta$ for which
\begin{eqnarray*}
v_J(i) & \geq & s_i + \beta c_{11}  \qquad\mbox{if }  i \in J,
\\
v_J(i) & \leq & s_i - \beta c_{11}  \qquad\mbox{if }  i \in I \setminus
J,
\end{eqnarray*}
and it is standard to verify that $(s_i)_{i \in I} \subset\beta
B_\infty^I$. It remains to show that $(X_i)_{i \in I}$ is
$c_{11}\beta$-shattered by $F$ itself. To that end, fix any $J
\subset I$, and let $f_J \in F$ be the function for which
\[
v_J= \sum_{\{i \dvtx |f_J(X_i)| \leq\beta\}} f_J(X_i)e_i + \sum_{\{i \dvtx |f_J(X_i)| > \beta\}} \beta\cdot\sgn{(f_J(X_i))} e_i.
\]
Observe that $(I \setminus J) \cap\{i\dvtx \sgn(v_J(i))>0\} \subset
\{i \dvtx |f_J(X_i)| \leq\beta\}$. Indeed, if there were some $i \in I
\setminus J$ for which $\sgn(v_J(i))>0$ and $|f_J(X_i)| > \beta$,
then on one hand, $v_J(i) =\beta$, but on the other, $v_J(i) \leq
s_i -\beta c_{11} \leq\beta(1-c_{11}) < \beta$, which is
impossible. In a similar fashion, $J \cap\{i\dvtx \sgn(v_J(i)) < 0\}
\subset\{i \dvtx |f_J(X_i)| \leq\beta\}$. Finally, fix $i \in J$. If
$f_J(X_i) \not= v_J(i)$ and $\sgn(v_J(i)) > 0$ then $f_J(X_i) \geq
v_J(i) \geq s_i + \beta c_{11}$. Otherwise, $\sgn(v_J(i)) < 0$,
implying that $v_J(i)=f_J(X_i)$. Hence, for every $i \in J$,
\[
f_J(X_i) \geq s_i + \beta c_{11},
\]
and by the same argument, for every $i \in I \setminus J$,
\[
f_J(X_i) \leq s_i - \beta c_{11}.
\]
Therefore, $\VC(F,\beta c_{11}) = \VC(F,c_{13}/\sqrt{k}) \geq c_{12}k$,
as claimed.
\end{pf*}

\section*{Acknowledgments}
The author would like to thank A. Libman, M. Kozdoda and the anonymous
referees for their careful reading of the manuscript and for many
valuable suggestions and comments.

The research leading to the results presented here has received funding
from the European Research\vadjust{\goodbreak} Council under the European Community's
Seventh Framework Programme (FP7/2007-2013)/ERC Grant
Agreement
[203134], from the Israel Science Foundation Grant 666/06 and
from the Australian Research Council Grant DP0986563.


%
\printaddresses

\end{document}